\documentclass[a4paper, 10pt, reqno]{amsart}

\usepackage{ amssymb, amsmath, amsthm}
\usepackage{graphicx, psfrag, setspace, subfigure, gensymb}
\usepackage{amsfonts}
\usepackage{bbm}
\usepackage{fix-cm}
\usepackage{comment}
\usepackage{tikz,tikz-cd}
\usetikzlibrary{matrix,arrows,decorations.pathmorphing,decorations.pathreplacing, decorations.markings}
\tikzset{commutative diagrams/diagrams={baseline=-2.5pt},commutative diagrams/arrow style=tikz}
\usepackage[colorlinks]{hyperref}
\usepackage{float}
\usepackage{setspace} 

\newcommand\Z{\mathbb Z}
\newcommand\C{\mathbb C}
\newcommand\R{\mathbb R}

\newcommand\A{\mathbb A}

\newcommand{\cA}{\mathcal{A}}

\newcommand{\cC}{\mathcal{C}}
\newcommand{\cD}{\mathcal{D}}
\newcommand{\cE}{\mathcal{E}}
\newcommand{\cF}{\mathcal F}
\newcommand{\cG}{\mathcal{G}}

\newcommand{\cI}{\mathcal{I}}

\newcommand{\cL}{\mathcal{L}}
\newcommand{\cN}{\mathcal{N}}
\newcommand{\cO}{\mathcal{O}}

\newcommand{\cZ}{\mathcal{Z}}

\newcommand{\set}[1]{\left\{{#1}\right\}}

\newcommand\isoto{\stackrel{\sim}{\To}}
\newcommand\id{\mathrm 1}

\newcommand\To{\longrightarrow}

\newcommand\Hom{\operatorname{Hom}}
\newcommand\End{\operatorname{End}}
\renewcommand\hom{\mathcal{H}om }

\renewcommand\P{\mathbb P}

\newcommand\Perf{\operatorname{Perf}}

\newcommand{\Crit}{\operatorname{Crit}}

\newcommand{\smat}[1]{\left(\begin{smallmatrix}#1\end{smallmatrix}\right)}

\newcommand\Spec{\operatorname{Spec}}

\newcommand{\al}[1]{\begin{align*}#1\end{align*}}

\newcommand{\beq}[1]{\begin{equation}\label{#1} }
\newcommand{\eeq}{\end{equation}}
\newcommand{\pgap}{\vspace{5pt}}

\theoremstyle{plain}
\newtheorem{prop}[equation]{Proposition}
\newtheorem{thm}[equation]{Theorem}
\newtheorem{lem}[equation]{Lemma}
\newtheorem{cor}[equation]{Corollary}

\theoremstyle{remark}
\newtheorem{rem}[equation]{Remark}

\theoremstyle{definition}
\newtheorem{defn}[equation]{Definition}

\newtheorem{eg}[equation]{Example}

\makeatletter \@addtoreset{equation}{section} \makeatother

\let\oldtocsection=\tocsection
\let\oldtocsubsection=\tocsubsection
\let\oldtocsubsubsection=\tocsubsubsection
\renewcommand{\tocsection}[3]{\hspace{0em}\oldtocsection{#1}{#2}{#3}}
\renewcommand{\tocsubsection}[3]{ \hspace{1em} \oldtocsubsection{#1}{\small{#2}}{\small{#3}} }
\renewcommand{\tocsubsubsection}[3]{\hspace{2em}\oldtocsubsubsection{#1}{\small{#2}}{\small{#3}}}

\newcommand{\marginparstretch}{0.6}
\let\oldmarginpar\marginpar
\renewcommand\marginpar[1]{\-\oldmarginpar[\framebox{\setstretch{\marginparstretch}\begin{minipage}{\marginparwidth}{\raggedleft\scriptsize #1}\end{minipage}}]{\framebox{\setstretch{\marginparstretch}\begin{minipage}{\marginparwidth}{\raggedright\scriptsize #1}\end{minipage}}}}

\newcommand{\aand}{\quad\quad\mbox{and}\quad\quad}

\newcommand\Fuk{\operatorname{Fuk}}
\tikzset{
       thin/.style=       {line width=2pt},
 every picture/.style={thin}
}

\newcommand\sslash{/\!/}
\newcommand\actson{\curvearrowright}

\begin{document}

\title{line fields on punctured surfaces and twisted derived categories}
\author{ Ed Segal}

\begin{abstract}
The Fukaya category of a punctured surface can be reconstructed from a pair-of-pants decomposition using a formal construction that attaches a category to a trivalent graph. We extend this formal construction to include a choice of line field on the surface, this requires a certain decoration on the graph. On the mirror side we show that this leads to a kind of twisted derived category which has not been widely studied.

Mirror symmetry predicts that our category should be an invariant of decorated graphs and we prove that this is indeed the case, using only B-model methods. We also give B-model proofs that a few different mirror constructions are equivalent.
\end{abstract}

\maketitle

\tableofcontents

\section{Introduction}

The simplest kind of symplectic manifold is a surface, and punctured surfaces are even simpler since they are exact. Because of this describing the Fukaya category of a punctured surface $\Sigma$ is a relatively tractable problem and there is an extensive literature on it \cite[...]{AAEKO, Bocklandt, DK, HKK, Lee, LPz, LPnodes, PS,  STZ}. 

 This problem is essentially the same problem as finding a homological mirror to $\Sigma$.  To `describe' $Fuk(\Sigma)$ means to come up with some construction of an equivalent category, and these other constructions are typically based on algebra or algebraic geometry. For example we might find a variety $X$ and take its derived category of coherent sheaves, or a non-commutative algebra $A$ and its derived category of modules, or a Landau-Ginzburg model $(X,W)$ and its category of matrix factorizations. If our algebro-geometric category is equivalent to $Fuk(\Sigma)$ then we can declare that $X$ (or $A$, or $(X, W)$,...) is a `mirror' to $\Sigma$.  This can be a rather facile version of mirror symmetry, especially for the more abstract algebraic constructions.

Contained in the list of references above are a wide range of different mirror constructions. Moreover even within a given construction there may be choices to be made, such as a choice of birational model for $X$. All the resulting categories are equivalent to $Fuk(\Sigma)$ so they must be equivalent to each other. This fact is not usually obvious on the algebraic side so it can be an interesting exercise to prove that two constructions are indeed equivalent using only B-model techniques, without invoking mirror symmetry. This is one of the themes of this paper.
\pgap

A second and perhaps more significant theme is to further understand the role of line fields in these mirror constructions. A line field is an extra piece of data on $\Sigma$, it's a section
 $$\eta \in \Gamma(\P T_\Sigma)$$
of the projectivized tangent bundle  (see Section \ref{sec.linefields}). For each choice of $\eta$ there is a $\Z$-graded version of the Fukaya category. Without a line field the Fukaya category is only $\Z_2$-graded.

  On the mirror side there is a well-known analogue of this for Landau-Ginzburg models: the category of matrix factorizations is initially only $\Z_2$-graded but we can get a $\Z$-graded category if we choose an appropriate $\C^*$ action on $X$. Such an action is called an \emph{R-charge} (Section \ref{sec.Rcharge}), locally it's a choice of grading on the ring $\cO_X$. 
The categories $D^b(X)$ or $D^b(A)$ are already $\Z$-graded but this corresponds to the trivial $\C^*$ action on $X$ or trivial grading on $A$;  other choices of R-charge are possible and produce different $\Z$-graded categories.

The main discovery of this paper is that R-charges are in general \emph{not sufficient} to mirror all the possible line fields on $\Sigma$. The remaining line fields are accounted for by the group $H^1_{\acute{e}t}(X, \Z)$. There is a simple but little-studied way to twist the derived category of $X$ by such an class (see Section \ref{sec.H1twists}) and we find that a general line field is mirror to some choice of R-charge together with some $H^1$-twist.
\pgap

More expert readers may wish to skip straight to Section \ref{sec.3pt} where we illustrate most of the contents of this paper using the example of a three-punctured torus. 

\subsection{Pair-of-pants decompositions}

There are many approaches to understanding $\Fuk(\Sigma)$, but for this paper the most important route is via \emph{pair-of-pants decompositions}. Any surface can be decomposed into a number of copies of the pair-of-pants (the three-punctured sphere) glued together along cylinders. Such a decomposition can be recorded as a trivalent graph $\Gamma$, with a vertex for each pair-of-pants, an internal edge for each cylinder, and an external edge for each puncture. Obviously $\Gamma$ records the genus and number of punctures so we can use it to reconstruct $\Sigma$ up to homeomorphism. Formally we replace each vertex with a pair-of-pants and each edge with a cylinder, add the inclusion maps, then we have diagram of spaces and the colimit is $\Sigma$.

Lee \cite{Lee} and Pascaleff-Sibilla \cite{PS} have independently proven that the Fukaya category of $\Sigma$ can be reconstructed in a similar way.\footnote{In this paper when we say `Fukaya category' we will always mean the \emph{wrapped} Fukaya category.}  The Fukaya category of the pair-of-pants admits three `restriction' functors to the Fukaya category of the cylinder, one for each leg. If we replace the vertices and edges of our graph by these categories and add the functors then we get a diagram of categories, and the limit of this diagram is $\Fuk(\Sigma)$.

This result provides one formal approach to the mirror of $\Sigma$. As input we need to know two basic instances of mirror symmetry:
\begin{enumerate} \item The Fukaya category of the cylinder is equivalent to the derived category of the punctured affine line $\A^1 \setminus 0$. 

\item The Fukaya category of the pair-of-pants is equivalent to the category of matrix factorizations of the polynomial $xyz \in \C[x,y,z]$. We denote this category by $D^b(\A^3, xyz)$. 
\end{enumerate}
 The three functors
\beq{eq.functorsintro}\begin{tikzcd}D^b(\A^3, xyz) \arrow[yshift=1ex]{r}\arrow[yshift=-1ex]{r}\arrow{r}\ar[r]& D^b(\A^1 \setminus 0)\end{tikzcd}\eeq
on the algebraic side come from Kn\"orrer periodicity (see Section \ref{sec.POPmirror}). Formally we can produce a mirror to $\Sigma$ simply by replacing each category in our diagram by its mirror category and then taking the limit. 

 In the abstract this procedure is vacuous, but in many cases we can give a worthwhile geometric description of this mirror limit category. For example, the oldest and most physically-meaningful version of mirror symmetry in this context is due to Hori and Vafa; for some $\Sigma$ the mirror is a toric Calabi-Yau threefold $X$ equipped with a certain canonical superpotential $W$. The dual to the toric fan is a planar trivalent graph whose vertices are the toric charts $\A^3 \subset X$.  In each chart the superpotential becomes $xyz$, so we get a diagram of categories just as above, and since matrix factorizations satisfy descent the limit category is $D^b(X, W)$.  This is how Lee and Pascaleff-Sibilla prove homological mirror symmetry for the Hori-Vafa mirrors.

\subsection{IH moves}\label{sec.IHmoves}

Given a surface $\Sigma$ and a pair-of-pants decomposition, we can construct a formal mirror category $\cD_\Gamma$ by taking the limit over the corresponding diagram. Since $\cD_\Gamma \cong Fuk(\Sigma)$, and the Fukaya category doesn't depend on the chosen decomposition, the category $\cD_\Gamma$ must actually only depend on the genus and number of external edges of $\Gamma$. This fact is not immediately obvious on the mirror side but it is possible to give an algebro-geometric proof, as we'll now sketch (for more details see Section \ref{sec.Z2invariance}). 
\pgap

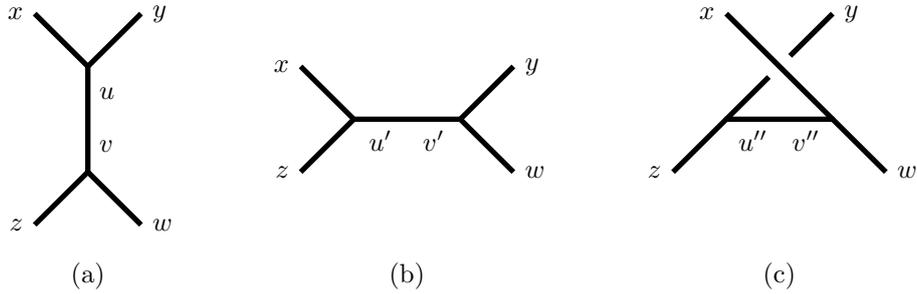
\begin{figure}
\begin{tikzpicture}[scale =.7]
\draw (0,0) --(1,1)--(1,3)--(0,4);
\draw (1,3)--(2, 4);
\draw (1,1)--(2,0);
\node [left] at (0,4) {$x$};
\node [right] at (2,4) {$y$};
\node [right] at  (1, 2.5) {$u$};
\node [right] at  (1, 1.5) {$v$};
\node [left] at  (0, 0) {$z$};
\node [right] at  (2,0) {$w$};

\node [below] at  (1,-.5) {(a)};

\begin{scope}[shift={(0,0)}]
\draw (5,1) --(6,2)--(8,2)--(9,1);
\draw (6,2)--(5, 3);
\draw (8,2)--(9,3);
\node [left] at (5,3) {$x$};
\node [right] at (9,3) {$y$};
\node [below] at  (6.5, 2) {$u'$};
\node [below] at  (7.5, 2) {$v'$};
\node [left] at  (5, 1) {$z$};
\node [right] at  (9,1) {$w$};
 
\node [below] at  (7,-.5) {(b)};
  \end{scope}

\begin{scope}[shift={(7,0)}]
\draw (5,1)--(6,2)--(8,2)--(9,1);
\draw (8,2)--(6, 4);
\draw (8,4)--(7.2,3.2);
\draw (6.8,2.8)--(6,2);
\node [left] at (6,4) {$x$};
\node [right] at (8,4) {$y$};
\node [below] at  (6.5, 2) {$u''$};
\node [below] at  (7.5, 2) {$v''$};
\node [left] at  (5, 1) {$z$};
\node [right] at  (9,1) {$w$};

\node [below] at  (7,-.5) {(c)};
 
  \end{scope}
\end{tikzpicture}
\caption{IH moves.}
\label{fig.IHmoves}
\end{figure}

For a given surface there are many possible pair-of-pants decompositions. For example a 4-punctured sphere admits exactly three decompositions, the corresponding trivalent graphs are in Figure \ref{fig.IHmoves}. This is related to an operation in graph theory called an \emph{IH move}: given a trivalent graph $\Gamma$ and a choice of internal edge we can get a new trivalent graph by contracting the edge and then expanding it again. In fact there are two possible new graphs since in Figure \ref{fig.IHmoves} we can replace (a) by either (b) or (c). It is an elementary fact that any two trivalent graphs with the same genus and same set of external edges are connected by a sequence of IH moves, so to prove that $\cD_\Gamma$ is a topological invariant of $\Gamma$ it's sufficient to prove that it doesn't change under IH moves. In fact it's enough to prove this when $\Gamma$ is a graph from Figure \ref{fig.IHmoves} since when we do an IH move only this part of the graph changes.

The 4-punctured sphere has a Hori-Vafa mirror, it's the small resolution $X$ of the 3-fold ODP singularity, equipped with a certain canonical superpotential $W$. This variety is covered by two toric charts and its fan is dual to the trivalent graph in Figure \ref{fig.IHmoves}(a). Passing to Figure \ref{fig.IHmoves}(b) corresponds to performing a standard flop, replacing $X$ with the other small resolution $X'$. A famous theorem of Bondal and Orlov says that derived categories do not change under flops, and the same is known to be true for categories of matrix factorizations, so $D^b(X, W)\cong D^b(X', W')$. This proves the invariance of $\cD_\Gamma$. 
\pgap

If $\Sigma$ has one Hori-Vafa mirror then it will have several, all related to each other by standard flops, and hence all having equivalent categories of matrix factorizations. The fact that $\cD_\Gamma$ is invariant under IH moves is a generalization of this observation. Note that flopping a toric CY 3-fold corresponds to an IH move which preserves the planar structure, which a general IH move does not.\footnote{See the recent \cite{PS2} for some discussion of some non-planar cases using $D_{sg}$ instead of $MF$.}

\subsection{Nodal curves}

Another possible mirror to the pair-of-pants is the node:
$$ Z = \{xy=0\} \subset \A^2 $$
This description fits nicely with an SYZ-picture: if we view $\Sigma$ as a cylinder with one extra puncture then it has the structure of an $S^1$-fibration with one degenerate fibre (the one that hits the puncture). On the mirror side the node is an $S^1$-fibration with one degenerate fibre at the origin. 
\pgap

Algebraically the equivalence between $D^b(\A^3, xyz)$ and $D^b(Z)$ is another instance of Kn\"orrer periodicity.\footnote{This statement should be treated with some caution, see Remark \ref{rem.KPfails}.}  In certain examples it allows us to give a different kind of geometric description of the formal mirror category $\cD_\Gamma$. 

There are two obvious functors from $D^b(Z)$ to $D^b(\A^1\setminus 0)$ given by simply restricting to the two branches of the node - these correspond to two of the three functors \eqref{eq.functorsintro}, the third functor is less geometric. So for some $\Gamma$ we can simply glue copies of $Z$ together, along copies of $\A^1\setminus 0$, to give us a nodal curve $C$. Then we can conclude that $D^b(C) \cong \cD_\Gamma$ is a mirror to the corresponding surface.

 This gluing is only possible if $\Gamma$ has the property that every vertex meets at least one external edge, since at any vertex we only have two geometric edges that we can glue along. Such a $\Gamma$ can only have genus zero or one. The resulting curve $C$ is either a ring of rational curves, or a chain of rational curves with an $\A^1$ at both ends. These one-dimensional mirrors were described by Sibilla-Treumann-Zaslow \cite{STZ} who called them \emph{balloon rings} and \emph{balloon chains} (see Section \ref{sec.STZ}). 
\pgap

In fact \cite{STZ} also allowed some orbifold structure at the nodes and this was extended by Lekili-Polishchuk \cite{LPnodes} to a broader class of orbifold structure, resulting in a mirror construction that works for a surface of arbitrary genus. In Section \ref{sec.orbifoldnodes} we give the B-model argument (\emph{i.e.} without invoking mirror symmetry) that connects these orbifold nodal curves to a Hori-Vafa model. It comes from Kn\"orrer periodicity combined with the invariance of $
D^b(X,W)$ under `orbifold flops'.

\subsection{Our results}

A lot of this paper is concerned with upgrading the discussion of the previous three sections to include line fields on $\Sigma$, and hence $\Z$-graded categories on the mirror side.

A line field on the pair-of-pants is determined by its three winding numbers $\alpha_x, \alpha_y, \alpha_z$ around the punctures, and these must satisfy
$$\alpha_x+\alpha_y+\alpha_z =2 $$
This data is mirror to an R-charge (a grading) on $\C[x,y,z]$ such that $\deg(xyz)=2$. If we have a line field $\eta$ on a general surface $\Sigma$, and we choose a pair-of-pants decomposition, then we should record these integers at each vertex of our trivalent graph. However this is not enough information to reconstruct $\eta$, since if we take surfaces with line fields and glue them together then the line field on the glued surface is not unique (see Section \ref{sec.linefields}).

 On the mirror side this corresponds to the fact that the functors \eqref{eq.functorsintro} are not canonical, they are only defined up to shifts (see Section \ref{sec.POPmirror}). This means that the limit category is not well-defined without further information.
\pgap

In Section \ref{sec.categoriesfromgraphs} we formulate a combinatorial definition of a \emph{decoration} on a trivalent graph $\Gamma$. This is exactly the data needed to build a line field $\eta$ on the corresponding surface $\Sigma$ (up to certain Dehn twists), or on the mirror side to define a $\Z$-graded limit category $\cC_\Gamma$. 
\pgap

It is almost certainly true that our category $\cC_\Gamma$ is equivalent to the $\Z$-graded Fukaya category of $(\Sigma, \eta)$, and presumably this could be proven by reworking \cite{Lee} or \cite{PS} with line fields in place. We have not done this. Instead, we have proven various results about these categories $\cC_\Gamma$, purely on the B-side.

\begin{itemize}\setlength{\itemsep}{5pt}

\item In Section \ref{sec.STZ} we consider the mirror nodal curves $C$ that arise from certain pair-of-pants decompositions in genus zero or one. In genus zero we see that any line field on the surface is mirror to some  R-charge $\alpha:\C^*\actson C$. In genus one we see that must also include a class $\beta\in H^1_{\acute{e}t}(C, \Z)$, and then we have an equivalence
$$\cC_\Gamma \cong D^b_\alpha(C, \beta) $$
to a twisted version of the derived category of $C$ (see the following section). 

\item In Section \ref{sec.HoriVafa} we show that certain planar graphs with a certain kind of decoration lead to Hori-Vafa mirrors; \emph{i.e.}~$\cC_\Gamma$ is equivalent to the category of matrix factorizations $D^b_\alpha(X,W)$ for some R-charge $\alpha: \C^*\actson X$.

 Some other choices of decoration appear to produce a version of this category twisted by a class $\beta\in H^1_{\acute{e}t}(\Crit(W), \Z)$. 

\item In Section \ref{sec.nodalelliptic} we consider the special case of a 1-punctured torus $\Sigma$, whose mirror is an elliptic curve $C$ with one node. Mirror symmetry predicts that we should have an equivalence
\beq{eq.introFM}D^b_\alpha(C, \beta) \isoto D^b_\beta(C, -\alpha)\eeq
mirror to the obvious element of the mapping class group $MCG(\Sigma)$. We prove that the usual Fourier-Mukai transform from $C$ to its Jacobian produces this equivalence.

\item In Section \ref{sec.invariance} we address the question of the invariance of $\cC_\Gamma$. The simple claim that the category should depend only on the topology of the graph becomes more complex in the $\Z$-graded case, because the Fukaya category depends on $\eta$ as well as $\Sigma$.

The orbits of line fields under the action of $MCG(\Sigma)$ have been classified in terms of certain invariants \cite{LPinvs}. In genus zero or one these are rather simple but in higher genus they involve the Arf invariant. We show how to combinatorially extract these same invariants from a decorated graph, then we prove the mirror symmetry prediction that $\cC_\Gamma$ depends only on the topology of $\Gamma$ plus these invariants of the decoration. 

The outline of the proof follows the sketch in Section \ref{sec.IHmoves} but tracking the behaviour of decorations under IH moves turns out to be rather delicate. It also uses our equivalence \eqref{eq.introFM}.

\end{itemize}

\subsection{Twisted derived categories}\label{sec.H1twists}

A `twisted derived category' of a variety $X$ is a family of categories on $X$ which is locally equivalent to the usual derived category $D^b(X)$. Azumaya algebras provide a classical example - if $\cA$ is a sheaf of algebras on $X$ which is locally a matrix algebra then locally $\cA$ is Morita equivalent to $\cO_X$, but perhaps not globally. To\"en \cite{Toen} has developed a more general theory of derived Azumaya algebras and shown they are classified by the group:
$$H^1_{\acute{e}t}(X, \Z) \times H^2_{\acute{e}t}(X, \cO^*)$$
The intuitive explanation for this group is very simple and closely analogous to the fact that $H^1(X, \cO^*)$ classifies line-bundles. At each point in $X$ our family of categories is $D^b(\C)$, so the `transition functions' of the family take values in the automorphism group of $D^b(\C)$. Since this is a category it has two kinds of automorphisms:
\begin{itemize}\item 1-automorphisms, \emph{i.e.}~invertible endofunctors. These are just the shift functors so form the group $\Z$.
\item 2-automorphisms, \emph{i.e.}~automorphisms of the identity functor. This is the center of the category $D^b(\C)$ and equals $\C^*$. 
\end{itemize}
These two kinds of automorphims explain the two factors in the group above and give us two possible kinds of twists.

The words `twisted derived category' are typically used to refer to twists of the second kind. These are gerbes, and include classical Azumaya algebras. However it is twists of the first kind that are relevant for this paper. With this kind of twisted derived category the objects are locally given by chain-complexes of coherent sheaves, but the transition functions can (and must) include shifts. 

This kind of twist do not seem to be widely known, in fact we don't know any previous references where they appear `in nature' or where examples have been studied. One reason for this may be that taking \' etale cohomology with integer coefficients is rather peculiar and often gives the zero group. 

These twists become more interesting when combined with a non-trivial R-charge since the the transition functions of (for example) an ordinary vector bundle might involve variables of non-zero cohomological degree. Such a vector bundle then becomes an object of the twisted derived category. See Remark \ref{rem.bmodalpha} or Section \ref{sec.nodalelliptic} for some explicit examples of this form.

\subsection{An example: the three-punctured torus}\label{sec.3pt}

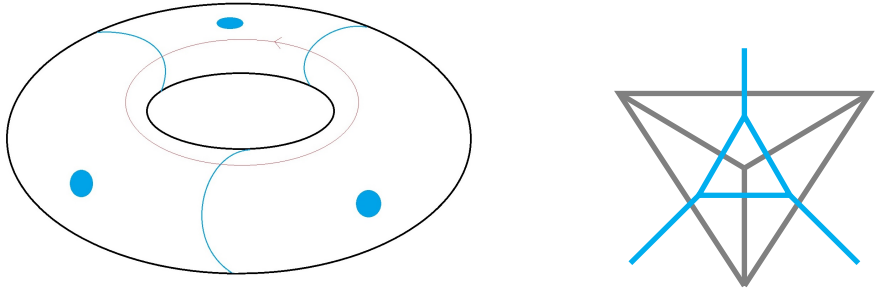
\begin{figure}[!tbp]
  \centering
  \begin{minipage}[b]{0.5\textwidth}
      \includegraphics[width=\textwidth]{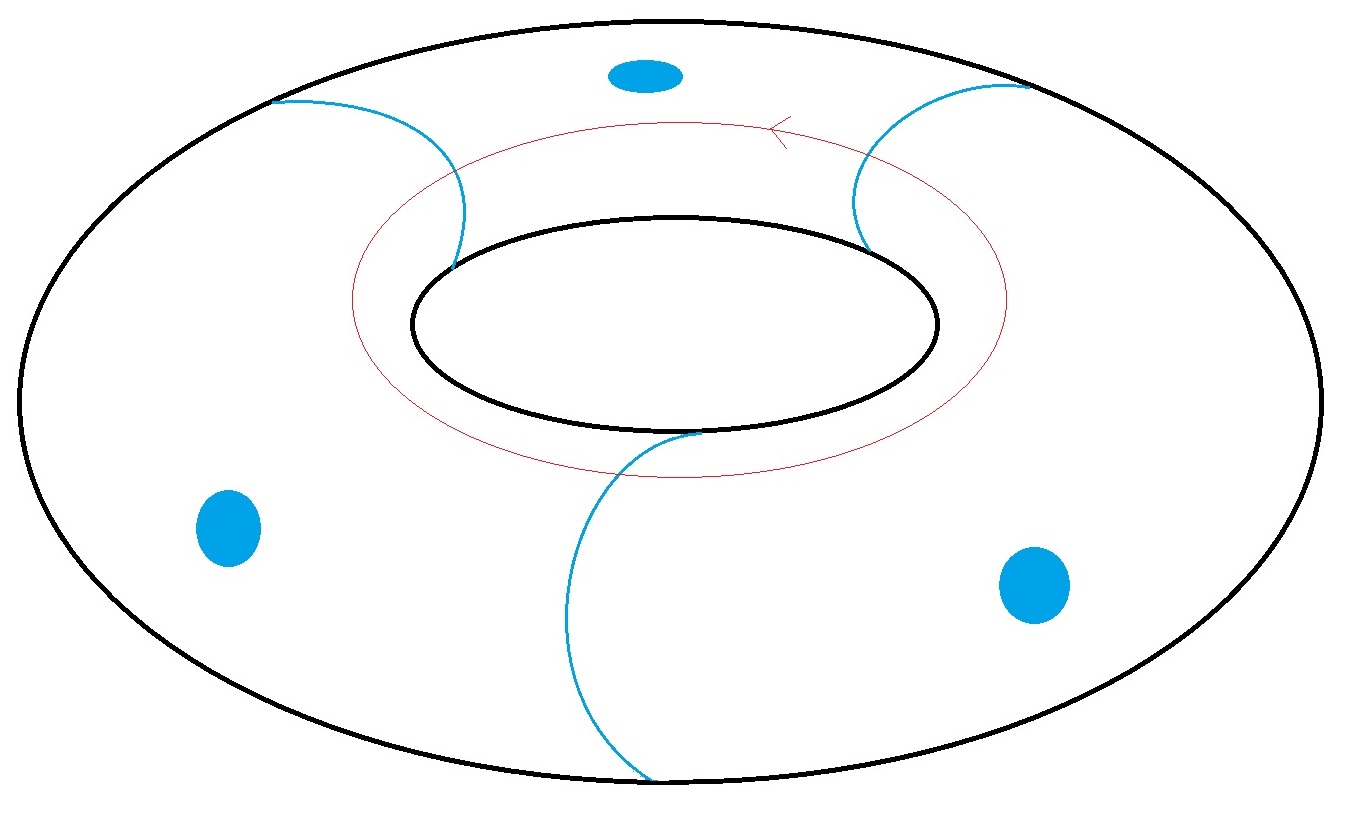}
  \end{minipage}
  \hfill
  \begin{minipage}[b]{0.35\textwidth}
    \begin{tikzpicture}[scale =.3]
	\draw[gray, thin] (2, 1.2)--(2, -4);
	\draw[gray, thin] (2,-4)--(7.5, 4.5)--(-3.5, 4.5)--(2,-4);
	\draw[gray] (2,1.2)--(7.5, 4.5);
	\draw[gray] (2,1.2)--(-3.5, 4.5);
	\draw [line width=.7mm, cyan] (-3,-3)--(0,0) --(4,0)--(2,3.5)--(0,0);
	\draw [line width=.7mm, cyan](4,0)--(7,-3);
	\draw [line width=.7mm, cyan](2,3.5)--(2,6.5);
   \end{tikzpicture}
  \end{minipage}
\caption{(L) A three-punctured torus with a given pair-of-pants decomposition. (R) The associated trivalent graph \emph{(in blue)} and the toric polytope of the mirror 3-fold \emph{(in grey)}. }
\label{fig.3punctured}
\end{figure}

Let $\Sigma$ be a torus with three punctures, as drawn on the left-hand-side of Figure \ref{fig.3punctured}. Cutting along the blue circles produces a pair-of-pants decomposition with three pieces. The associated trivalent graph $\Gamma$ is a triangle with an external edge emanating from each vertex, as shown on the right-hand-side. 
\pgap

The dual of this graph $\Gamma$ is a triangulation of a triangle with one interior point. We recognize this as the toric polytope for the Calabi-Yau threefold $K_{\P^2}$. 

On $K_{\P^2}$ we can put a superpotential $W=xyzp$ where $x,y,z$ are the projective co-ordinates on $\P^2$ and $p$ is the fibre co-ordinate. We have three toric charts, each of which is isomorphic to $\A^3$, and on each one $W$ is the product of the three co-ordinates. So $(K_{\P^2}, W)$ is three copies of the mirror to the pair-of-pants, glued together following $\Gamma$. Hence $\Fuk(\Sigma)\cong \cD_\Gamma$ is the category of matrix factorizations $D^b(K_{\P^2}, W)$. This is the Hori-Vafa mirror to $\Sigma$. 
\pgap

Alternatively, we take the mirror to each pair-of-pants to be a node. Gluing the branches of the nodes together we obtain a genus one curve $C$ which is a triangle of three rational curves meeting at nodes. It can be realized $C = V(xyz)\subset \P^2$. 

This is the `balloon ring' mirror to $\Sigma$. On the B-side the equivalence between $C$ and the Hori-Vafa mirror is an instance of global Kn\"orrer periodicity. 
\pgap

A third possible mirror is provided by \cite{LPnodes}, it is a node $Z = V(xy)\subset \A^2$ orbifolded by $\Z_3$ acting with weight 1 on each co-ordinate. This mirror doesn't use any pair-of-pants decomposition. To get from here to the Hori-Vafa mirror we first use equivariant Kn\"orrer periodicity to pass to the orbifold $[\A^3 / \Z_3] $ with the superpotential $W'=xyz$, where the $\Z_3$ also acts with weight 1 on $z$. Then we perform an `orbifold flop' to pass to the Hori-Vafa mirror $(K_{\P^2}, W)$ and this preserves the categories of matrix factorizations. 

In total we have four possible mirrors, and we have algebro-geometric equivalences connecting all four:
\beq{eq.4differentmirrors}D^b(C) \;\stackrel{Kn\ddot{o}}{\cong}\; D^b(K_{\P^2}, xyzp) \;\stackrel{flop}{\cong}\; D^b([\A^3/\Z_3], xyz)\; \stackrel{Kn\ddot{o}}{\cong} \;D^b([Z/\Z_3]) \eeq

\begin{rem}Typically a toric CY fold will have many birational models, some of which are orbifolds, all of which give equivalent categories of matrix factorizations. But it's only in special cases that we can apply Kn\"orrer periodicity globally to get a 1-dimensional mirror.
\end{rem}

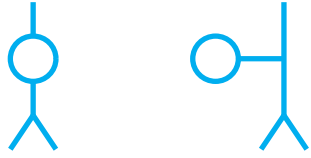
\begin{figure}[!tbp]
\centering
    \begin{tikzpicture}[scale =.3]
	\draw [line width=.7mm, cyan] (-1,1)--(-1,2.5);
       \draw [line width=.7mm, cyan] (0,0) arc [radius=1, start angle=0, end angle= 360];
      \draw [line width=.7mm, cyan] (-1, -1)--(-1, -2.5);
      \draw [line width=.7mm, cyan] (-1, -2.5)--(-2, -4);
      \draw [line width=.7mm, cyan] (-1, -2.5)--(0, -4);
 
	\draw [line width=.7mm, cyan] (10,-2.5)--(10,2.5);
       \draw [line width=.7mm, cyan] (8,0) arc [radius=1, start angle=0, end angle= 360];
      \draw [line width=.7mm, cyan] (10, 0)--(8, 0);
      \draw [line width=.7mm, cyan] (10, -2.5)--(9, -4);
      \draw [line width=.7mm, cyan] (10, -2.5)--(11, -4);
    \end{tikzpicture}
\caption{Other pair-of-pants decompositions for the three-punctured torus.}
\label{fig.genus1-3ext}
\end{figure}

There are two more trivalent graphs of genus one with three external edges, they're pictured in Figure \ref{fig.genus1-3ext} and correspond to other pair-of-pants decompositions of $\Sigma$. We can produce the mirror as an abstract category by gluing together three copies of $D^b(\A^3, xyz)$ following either of these graphs, but we don't know a global geometric interpretation in either case.
\pgap

Now we discuss the role of line fields in this example. A line field on $\Sigma$ is characterised by four integers $\alpha_1, \alpha_2, \alpha_3, \beta \in \Z$. The first three are the winding numbers along the blue loops in Figure \ref{fig.3punctured}(L)  - the loops we cut to produce the pair-of-pants decomposition - and the fourth is the winding number along the the red loop. 
 
In a single pair-of-pants we see two of the $\alpha_i$ as winding numbers around the two `cuffs'. On the mirror side these numbers give an R-charge (a grading) on the corresponding copy of $\A^3$. 

Globally, if we take the nodal curve $C$ as our mirror then we see that all three of the $\alpha_i$ can be accounted for by an R-charge $\alpha:\C^*\actson C$ since we can rotate the three rational curves with three different weights.  The fourth winding number is more interesting, it's mirror to a class:
$$\beta\in H^1_{\acute{e}t}(C, \Z) = \Z$$
 If $\beta \neq 0$ then the mirror to $(\Sigma, \eta)$ is a twisted derived cateory $D^b_\alpha(C, \beta)$. If we were to cut one of the blue loops then we'd get a 5-punctured sphere, encoded by a genus zero graph $\Delta$, and both the fourth winding number and the \'{e}tale cohomology group disappear. On both sides of the mirror $\beta$ is telling us how to glue the $\Z$-graded category $\cC_\Delta$ to itself. 
\pgap

If we take the Hori-Vafa mirror instead then we are more restricted about which line fields we can find mirrors to. Not all of the R-charges on $C$ extend to the three-fold $K_\P^2$; there is a rank 3 torus action but $W$ is only preserved by a rank 2 subtorus. In fact an R-charge on $(K_{\P^2}, W)$ can only be mirror to a line field satisfying $\alpha_1+\alpha_2+\alpha_3=0$. The fourth winding number $\beta$ is also harder to interpret in this model, but notice that the critical locus $\Crit(W)$ is the curve $C$ with three $\A^1$'s attached at the nodes, so:
$$H^1_{\acute{e}t}(\Crit(W), \Z)=H^1_{\acute{e}t}(C, \Z) = \Z$$
Since matrix factorizations are supported on $\Crit(W)$ it should be possible to twist the category by a class in $H^1_{\acute{e}t}(\Crit(W), \Z)$. 

If we consider the remaining two mirrors from \eqref{eq.4differentmirrors} then the same restriction on R-charge applies, and we are not certain how to interpret $\beta$.

\subsection{Desiderata}

There is another approach to mirror symmetry for punctured surfaces developed in \cite{Bocklandt, HKK, DK} where we cut $\Sigma$ into polygons instead of pairs-of-pants. This is a rather different procedure. We cut along arcs joining punctures, and each polygon is interpreted as a disc with stops on the boundary corresponding to the arcs where we cut. The Fukaya category of this disc-with-$n$-stops is the derived category of the $A_n$ quiver algebra. We glue the discs along rectangles, \emph{i.e.}~discs with 2 stops, but now the functors go in the opposite direction to the case of pair-of-pants decompositions and we recover $\Fuk(\Sigma)$ as the colimit (rather than the limit) of the diagram. This distinction makes it tricky to unify the two approaches.

For example the pair-of-pants can be cut into two hexagons, which we interpret as discs with 3 stops. It follows that $D^b(\A^3, xyz)$ can be built from two copies of $D^b(A_2)$  - which is equivalent to $D^b([\A^1/\Z_3], x^3)$ - glued together along three exceptional objects. We do not have a nice interpretation of this fact in algebraic geometry. 
\pgap

In Bocklandt's version \cite{Bocklandt} he interprets the cut surface as a dimer model and constructs the `mirror' as a dual dimer model, giving a category of matrix factorizations over a non-commutative algebra. In some cases this fits into the well-studied story of non-commutative crepant resolutions for toric 3-fold singularities. For example if we apply this to the example of the 3-punctured torus from the previous section then we get a non-commutative algebra $A$, but it turns out to be the twisted group algebra $\C[\Z_3]\rtimes \C[x,y,z]$. So this is equivalent in a rather trivial way to the third mirror in our previous list \eqref{eq.4differentmirrors}. In other examples Bocklandt's mirrors are genuinely different. 

Some line fields on $\Sigma$ will be mirror to gradings on the algebra, but not all. In the example of the 3-punctured torus there is a rank 5 set of gradings on $A$ but a rank 3 subset come from inner derivations, so we hit exactly the same restriction on the $\alpha_i$'s as before. We also don't know how to interpret $\beta$ here.
\pgap

Another direction is to consider compact surfaces.  Pascaleff-Sibilla recently \cite{PS2} extended their approach to this case. Since our primary interest is in the role of line fields the only compact surface we would care about is the torus, see Remark \ref{rem.compacttorus}.

\subsection{Acknowledgements}

It's a pleasure to thank Yank{\i} Lekili for his huge influence on this paper. It's basically a write-up of many long discussions between the two of us.

This project has received funding from the European Research Council (ERC) under the European Union Horizon 2020 research and innovation programme (grant agreement No.725010).

\section{Background}

\subsection{Gluing line fields}\label{sec.linefields}

Let $\Sigma$ be a symplectic surface with punctures, \emph{i.e.} with at least one boundary component. A \textit{line field} on $\Sigma$ is a section $\eta$ of the projectivised tangent bundle $\P T_\Sigma$.  We will need a few facts about line fields, we refer to \cite[Sec.~1]{LPinvs} for background and further references.
\pgap

We only care about line fields up to homotopy, and the homotopy classes form a torsor for
$$H^1(\Sigma,\Z)= \Z^{2g+n-1} $$
where $g$ and $n$ are the genus and number of punctures of $\Sigma$. To any embedded curve $\sigma$ in $\Sigma$ we can assign a winding number $w_\eta(\sigma)\in \Z$ by comparing $\eta$ to the tangent line field along $\sigma$. In particular each boundary component has an associated winding number. By the Poincar\'e-Hopf theorem the sum of these boundary winding numbers must be $2\chi(\Sigma)$. 
\pgap

If choose a line field $\eta$ on $\Sigma$ then we can define a $\Z$-graded version $\Fuk_\eta(\Sigma)$ of the wrapped Fukaya category. 
 
The usual wrapped Fukaya category is $\Z_2$-graded, and its objects are either compact curves or arcs which begin and end at the boundary. In $\Fuk_\eta(\Sigma)$ we only consider \emph{gradeable} Lagrangians, meaning either arcs, or curves such that $w_\eta(\sigma)=0$. This restriction means that $\Fuk_\eta(\Sigma)$ is `smaller' than the usual $\Z_2$-graded version.

\begin{eg}\label{eg.cylinder1}
Let $\Sigma$ be a cylinder, \emph{i.e.} a twice-punctured sphere. A line field on $\Sigma$ is classified by its winding number $a\in \Z$ around the the cylinder.  For $a\neq 0$ there are no compact gradeable Lagrangians. 
\end{eg}

A \emph{grading} on a gradeable Lagrangian is a choice (up to homotopy) of a homotopy between the tangent line field and $\eta$. The set of such homotopies forms a $\Z$-torsor and the shift operator on $\Fuk_\eta(\Sigma)$ acts by modifying the grading by 1 for every object.
\pgap

Now suppose we cut $\Sigma$ along a separating embedded curve $\sigma$, splitting it into surfaces $\Sigma_1$ and $\Sigma_2$, each of which has $\sigma$ as a boundary component.  If we have line fields $\eta_1$ and $\eta_2$ on the two pieces we can attempt to glue them to form a line field $\eta$ on $\Sigma$. Performing this gluing requires a choice of homotopy between $\eta_1|_{\sigma}$ and $\eta_2|_{\sigma}$, so it is only possible if the winding numbers $w_{\eta_1}(\sigma)$ and $w_{\eta_2}(\sigma)$ agree.  If $\sigma$ is non-separating essentially the same comments apply: cutting along $\sigma$ produces a single surface $\Sigma'$ with two extra boundary components $\sigma_1$ and $\sigma_2$, and a line field $\eta'$ on $\Sigma'$ can be glued to a line field $\eta$ on $\Sigma$ iff $w_{\eta'}(\sigma_1)=w_{\eta'}(\sigma_2)$. 
 
For non-separating $\sigma$ the glued line field $\eta$ is not unique. If $\eta'|_{\sigma_1}$ and $\eta'|_{\sigma_2}$ are homotopic then the space of homotopies between them has countably-many connected components, each of which produces a different homotopy class for $\eta$. For example suppose $\Sigma$ is a torus and we cut it to produce a cylinder $\Sigma'$; a line field on $\Sigma$ depends on two integers $(\alpha, \beta)$ but the cylinder only sees $\alpha$.

\begin{rem} 
Let $\cN_i$ be a tubular neighbourhood of one of the two boundary components $\sigma_i$ in $\Sigma'$, so each $\cN_i$ is a cylinder. There is a functor from $\Fuk_\eta(\Sigma')$ to $\Fuk_\eta(\cN_i)$ which just maps any Lagrangian to its intersection with $\cN_i$. 

We can recover our uncut surface $\Sigma$ by making an identification $\cN_1\simeq \cN_2\simeq \cN$ between these two boundary cylinders. Then we can try to recover $\Fuk_\eta(\Sigma)$ as the equalizer of the diagram:
$$\Fuk_{\eta'}(\Sigma') \rightrightarrows \Fuk_{\eta'}(\cN) $$
But there is an ambiguity here because the equivalence $\Fuk_{\eta'}(\cN_1)\cong \Fuk_{\eta'}(\cN_2)$ is only defined up to shifts. This is exactly the same ambiguity that arises in gluing the line field. 
\end{rem}

We can fix this ambiguity in the gluing by `framing' our surfaces as follows:
\begin{itemize} \item On each component of $\partial\Sigma$ we choose a marked point $p_i$.\footnote{We should emphasise that these marked points are \emph{not} stops for the Fukaya category - they have no effect on the Fukaya category.}
\item A line field $\eta$ on a $\Sigma$ must obey the condition that $\eta$ is not tangent to $\partial \Sigma$ at any marked point $p_i$. A homotopy of line fields must obey this condition at all times. 
\item When we cut a surface along a curve $\sigma$ we should also choose a marked point on $\sigma$ to become a marked point on the two new boundary components. 
\item When we glue surfaces along boundary components we we glue the marked points together.
\end{itemize}
With these conventions gluing of line fields becomes well-defined: if $w_{\eta'}(\sigma_1)=w_{\eta'}(\sigma_2)$ then there is a unique homotopy between $\eta'|_{\sigma_1}$ and $\eta'|_{\sigma_2}$ which is never tangent to the curve at the marked point.
\pgap

On a framed surface the set of line fields up-to-homotopy (with the conventions above) is a torsor for the relative cohomology group $H^1(\Sigma, \set{p_1,..., p_k})$.  We can assign winding numbers not just to closed curves in $\Sigma$, but also to arcs $\sigma$ which connect two of the marked points, provided that they meet $\partial\Sigma$ transversely. 
\pgap

An important example for us is when $\Sigma$ is the  pair-of-pants. An ordinary line field on $\Sigma$ is characterized by its winding numbers $\alpha_1, \alpha_2, \alpha_3\in \Z$ about the three boundary components, which must obey:
$$\alpha_1+\alpha_2 + \alpha_3=2\chi(\Sigma)=2$$
 If we now add marked points $p_1, p_2, p_3$ on each boundary component then it takes a further two integers to specify a line field, these can be detected as winding numbers along arcs connecting the $p_i$. 

More precisely, suppose $\sigma$ is an embedded arc from $p_1$ to $p_2$ which meets $\partial\Sigma$ transversely. A line field $\eta$ has a winding number $\widetilde{\gamma}_{12}\in \Z$ along this arc. But there are of course many such arcs -  we can act by a Dehn twist along the boundary component at either end and this will change the winding number by $\alpha_1$ or $\alpha_2$. Since all possible arcs (up to homotopy) are obtained this way the well-defined invariant is the residue:
$$\gamma_{12} = [\widetilde{\gamma}_{12}]\in \Z/(\alpha_1, \alpha_2)$$ 
There are six such invariants, one for each ordered pair of boundary components. Obviously $\gamma_{ij}=-\gamma_{ji}$ and there also the relation
$$\gamma_{12} + \gamma_{23} + \gamma_{31} = 1 \in \Z/(\alpha_1, \alpha_2, \alpha_3) $$
which can be seen by considering the case when three arcs bound a triangle.

\begin{rem}\label{rem.movep}
We can also ask what happens if we fix a line field and move the marked points. If $\eta$ has winding number $\alpha_1$ along the first boundary component $\sigma_1$, then the subset of $\sigma_1$ where $\eta$ is not tangent to the boundary must have at least $\alpha_1$ connected components. If we move $p_1$ through a point where $\eta$ crosses the boundary then $\gamma_{12}$ and $\gamma_{13}$ increase by 1. 
\end{rem}

We can now address the question of line fields on a pair-of-pants decomposition. Take a surface $\Sigma$ and decompose it into pairs-of-pants by cutting along curves $\sigma_1, ..., \sigma_k$. The decomposition can be recorded as a trivalent graph and this combinatorial object is enough to reconstruct $\Sigma$. 

Now suppose we also want to record the data of a line field $\eta$ on $\Sigma$. We have a winding number $\alpha_i$ along each curve $\sigma_i$ and we can label each internal edge of the graph with this integer. More precisely we should label half of the edge with $\alpha_i$, and the other half with $-\alpha_i$, to record the orientation correctly.

This data specifies a line field on each pair-of-pants but it is not enough to reconstruct $\eta$ since, as discussed above, the glued  line field is not unique. We must also choose a marked point $p_i$ on each $\sigma_i$, and record our $\gamma$ invariants for each pair-of-pants. 

If we label the half-edges of our graph with the $\alpha$'s and the vertices with the $\gamma$'s then we have enough information to reconstruct the line field $\eta$, up to the action of Dehn twists around the curves $\sigma_i$.  In Section \ref{sec.categoriesfromgraphs} we will codify this kind of labelled graph precisely.

\subsection{R-charge and matrix factorizations}\label{sec.Rcharge}

In this paper we will frequently deal with rings $S$ equipped with a grading. These gradings will be 'cohomological' gradings rather than `internal' gradings, \emph{i.e.}~we want to regard $S$ as a dga that happens to have zero differential. Then the obvious way to define the derived category $D(S)$ would be to take dg-modules, that is graded $S$-modules with a differential of  degree 1 (rather than chain-complexes of graded modules). 

The global analogue of a grading is a $\C^*$ action
$$\alpha: \C^* \actson X$$
on a scheme. The global analogue of a dg-module is an equivariant coherent sheaf $\cE$ on $X$  equipped with a differential 
$d: \cE \to \cE[1]$.  Here the operator $[n]$ means twisting by the $n$-th character of $\C^*$.  

We call the grading or the $\C^*$ action the \emph{R-charge}, it's the choice of vector R-charge in the N=(2,2) supersymmetric field theory with target $X$. 
\pgap

One can define a version of the derived category with this approach, let's call it $D^b_\alpha(X)^{obv}$. We're not going to use it because it has a serious technical defect, as we now explain. 

 Locally $\cO_X$ is a commutative graded ring. However, it is \emph{not} supercommutative unless the grading happens to be concentrated in even degrees. Globally this is the hypothesis
\beq{eq.SC} \alpha(-1) = \id_X \eeq
\emph{i.e.}~the R-charge is the square of some other $\C^*$ action. We do not want to make this hypothesis since on the mirror side the degrees correspond to (some of) the winding numbers of the line field $\eta$, and these are not necessarily even.\footnote{In fact all winding numbers are even iff $\eta$ lifts to a vector field.}

Without supercommutativity we do not have internal homs in our category; if $\cE$ and $\cF$ are objects of $D^b_\alpha(X)^{obv}$ then we have a sheaf $\hom(\cE, \cF)$ but, because of the signs, the differential is not linear over the structure sheaf.  This means that we cannot localize the category to open sets in $X$ (or at least not obviously), so we cannot rebuild the global category from the categories associated to an open cover. Since this is one of the key tools in this paper we will need to use a different definition. 
\pgap

This other definition starts with an idea which is common in the study of matrix factorizations. Classically the category of matrix factorizations is only $\Z_2$-graded (or $\Z_2$-periodic), but it is possible to define a $\Z$-graded category if we have some additional choice of grading or group action. 

\begin{defn}[e.g.~\cite{BFK}]\label{def.MF}
Fix a scheme $X$, a torus action $T\actson X$, a character $\chi$ of $T$, and a semi-invariant function:
$$ W \in \Gamma(\cO_X(\chi))$$
A \emph{matrix factorization} is a pair of $T$-equivariant coherent sheaves and maps
$$\cE_0 \stackrel{d_0}{\To} \cE_1 \stackrel{d_1}{\To} \cE_0(\chi)$$
such that $d_1d_0=W\id_{E_0}$ and $d_0d_1=W\id_{E_1}$. 
\end{defn}
$W$ is called the \emph{superpotential}. Matrix factorizations form a dg-category which we'll denote by:
$$D^b_T(X, W)$$
The definition of the morphisms takes some care and we won't recall the details since we won't need them.

\begin{rem}\label{rem.MFs}\,
\begin{enumerate}\setlength{\itemsep}{7pt}

\item In this category we take all morphisms and not just $T$-equivariant ones. Hence the hom-spaces are graded by the lattice
$$L = \frac{\Hom(T, \C^*)\times \Z }{ (\chi, -2)}$$
because the shift functor satisfies $[2] =\cO(\chi)$. 

\item One can replace $[X/T]$ and the line-bundle $\cO(\chi)$ with a more general stack equipped with a line-bundle, but this generality is sufficient for us.

\item If $W\neq 0$ then it determines $\chi$ but if $W=0$ then $\chi$ is extra data (even though we've suppressed it from the notation).

 We can recover the usual derived category $D^b(X)$ by choosing $T=\C^*$ with the trivial action on $X$, setting $W=0$, and setting $\chi=1$. Then $\cE_0$ and $\cE_1$ are just the even and odd summands of a chain-complex. 

\item Generalizing (3), if we have a torus action $T'\actson X$ we can recover the $T'$-equivariant derived category of $X$. We set $T=T'\times \C$, let the second factor act trivially, then set $W=0$ and $\chi=(0,1)$. In fact any $\chi$ of the form $(\chi', 1)$ will do; all such categories are equivalent and have canonically-isomorphic grading lattices.

\item In \cite{segal} there is a definition of a $\Z$-graded matrix factorization category from a scheme $X$ equipped with an R-charge $\alpha: \C^*\actson X$ satisfying \eqref{eq.SC}, plus a superpotential $W$ of weight 2. There a matrix factorization is a single equivariant sheaf $\cE$ equipped with a differential $d$ of weight 1, such that $d^2=W$. This generalizes the objects in $D^b_\alpha(X)^{obv}$ discussed above.

The condition \eqref{eq.SC} implies that such an $\cE$ splits into even and odd summands $\cE_0\oplus \cE_1$,  so we can recover this definition from Definition \ref{def.MF} by setting $T=\C^*$, acting with the square-root of $\alpha$, and setting $\chi=1$. 

\item Given two objects $\cE, \cF \in D^b_T(X,W)$ the `sheaf hom' between them give a matrix factorization $\hom(\cE, \cF)$ for $(X, T, \chi, 0)$.

\end{enumerate}

\end{rem}

We now have a category graded by some lattice $L$ as in Remark \ref{rem.MFs}(1) above.  Since we only wanted a $\Z$-graded category we now have to collapse this grading somehow.

\begin{defn}\label{defn.Rcharge} Given $(X, T, \chi, W)$ as in Definition \ref{def.MF} we define an \emph{R-charge} to be a homomorphism 
$$\alpha: \C^* \to T$$
such that $\chi(\alpha)=2$. 
\end{defn}

An R-charge gives us a homomorphism $(\alpha, 1): L \to \Z$ which we can use to collapse the grading on $D^b_T(X, W)$ to a $\Z$-grading. We define
$$D^b_\alpha(X, W)$$
to be the pre-triangulated closure of this $\Z$-graded dg-category.

\begin{rem}\hspace{1pt}\label{rem.Dbalpha}
 \begin{enumerate}\setlength{\itemsep}{7pt}

\item In the case of the usual derived category (Remark \ref{rem.MFs}(3)), or matrix factorizations as defined in Remark \ref{rem.MFs}(5),  $\alpha$ is forced to be the squaring map $\C^*\to \C^*$. It gives an isomorphism $L\isoto \Z$ so there is no change in the grading.

\item In the case of the equivariant derived category of some torus action $T'\actson X$, as in Remark \ref{rem.MFs}(4), an R-charge is equivalent to a homomorphism $\alpha': \C^* \to T'$. The grading lattice is
$$L = \Hom(T', \C^*) \times \left(\Z^2/(1, -2)\right) \cong \Hom(T', \C^*) \times \Z$$
and an R-charge gives a grading collapse $(\alpha', 1): L \to \Z$.

\item The step where we take the pre-triangulated closure is essential. Suppose $f,g:\cE\to \cF$ are two morphisms in $D^b_T(X, W)$ of different degrees, which become the same degree after we apply $\alpha$. Then in $D^b_\alpha(X,W)$ we want to be able to take the mapping cone on $f+g$, but there is no such object in $D^b_T(X,W)$. See Example \ref{eg.localisingA1} below. 

\item The approach we're using seems related to `universal gradings' on Fukaya categories, see \emph{e.g.}~\cite[Sec.~3]{Sheridan}.
\end{enumerate}
\end{rem}

\begin{eg}\label{eg.localisingA1}
Let $X=\A^1$ and let $T=(\C^*)^2$ acting on $X$ with weights $(1,0)$. Set $\chi=(0,1)$ and $W=0$. Then $D^b_T(X,W)$ is the usual equivariant derived category $D^b_{\C^*}(\A^1)$ (Remark \ref{rem.MFs}(4)). 

The $\C^*$ weights and the cohomological grading together grade this category by $L\cong\Z^2$. An R-charge is a homomorphism $\alpha=(a,2):\C^*\to (\C^*)^2$ for some $a\in \Z$. It's the same as a grading of the ring $\cO_X=\C[z]$ where we set $\deg(z)=a$. 

If $a$ is odd then this ring is not supercommutative, so in the obvious derived category $D^b_a(\A^1)^{obv}$ the hom spaces will not be linear over $\C[z]$. Then it's not immediately clear for example if there is a localization functor:
$$D^b_a(\A^1)^{obv} \To D^b_a(\A^1\setminus 0)^{obv}$$
But this localization functor certainly does exist for our $D^b_\alpha(\A^1)$, it's just induced from the localization on the equivariant derived category.

In the case $a=0$ we put our ring $\cO_X$ in degree zero and we end up with the usual derived category of $\A^1$. This case highlights the issue with mapping cones. The equivariant derived category of $\A^1$ only has torsion sheaves at the origin but $D^b(\A^1)$ has torsion sheaves at other points, for example point sheaves:
$$\big[ \cO \stackrel{z-\lambda}{\To} \cO \big] \isoto \cO_\lambda $$
We cannot form this mapping cone in $D^b_T(\A^1)$ so collapsing the grading gives a category which is not pretriangulated.
\end{eg}

\begin{rem}\label{rem.punctureA1} The punctured affine line $\A^1\setminus 0$ is the mirror to the cylinder. The choice of R-charge $\alpha=(a,2)$ is mirror to the choice of line field with winding number $a\in \Z$, see Example \ref{eg.cylinder1}. The torsion sheaves are mirror to compact Lagrangians, both exist only for $a=0$. 

For $a\neq 0$ the category $D^b_\alpha(\A^1\setminus 0)$ is extremely simple: all objects must be free modules over $\C[z,z^{-1}]$, and $z$ is an isomorphism between $\cO$ and $\cO[a]$, so this is the category of $\Z_a$-graded vector spaces. 
\end{rem}

\begin{rem} For the case $X=\A^1\setminus 0$ just discussed it's easy to show that $D^b_\alpha(X)$ agrees with the more straightforward definition  $D^b_\alpha(X)^{obv}$, for any R-charge $\alpha$. It seems plausible that this in fact holds for any smooth scheme $X$. It definitely does not hold for singular $X$, see Remark \ref{rem.KPfails} below.
\end{rem}

\begin{rem} As discussed in Section \ref{sec.H1twists} we will be interested in versions of our categories twisted by classes in $H^1_{\acute{e}t}(X, \Z)$. However we will only need these for the special case of equivariant derived categories, as in Remark \ref{rem.MFs}(4), and these are included in T\"oen's theory \cite{Toen} since he works over stacks. Presumably the theory could be generalized, see Remark \ref{rem.twistedCrit}.
\end{rem}

\subsection{Kn\"orrer periodicity and VGIT}

Many of our arguments are based on two kinds of important equivalences that exist between matrix factorization categories on different varieties.  The first of these is \emph{Kn\"orrer periodicity}. 
\pgap

Let $Y=\Spec S$ be a smooth variety with a torus action $T\actson Y$ and let $X=\mbox{Tot}\{\cL \stackrel{\pi}{\to} Y\}$ be the total space of a $T$-equivariant line-bundle. Let $f$ be a non-zero section of $\cL^\vee(\chi)$ for some character $\chi$, it induces a $\chi$-semi-invariant superpotential $W=fp$ on $X$ ($p$ is the fibre co-ordinate).

  Let $Z= (f)\subset X$ be the zero locus of $f$.  We equip it with the zero superpotential, viewed as a section of $\cO(\chi)$. 

Write $\Crit(W)$ for the critical locus of $W$. If $Z$ is smooth then $\Crit(W)=Z$, if not it includes the $p$ fibres at the singular points.

\begin{thm}\cite[...]{OrlovKP, Shipman, Isik, Hirano}\label{thm.KP} Assume that there are no $T$-invariant closed subvarieties in $\Crit(W)$ disjoint from the zero section. Then for any R-charge $\alpha$ we have an equivalence:
$$D^b_\alpha(X, W) \isoto D^b_\alpha(Z)$$
\end{thm}
Note that if $Z$ is smooth then the assumption holds automatically, but if $Z$ is singular it might not. See Section \ref{sec.POPmirror} for an example of this result and Remark \ref{rem.KPfails} to see why the assumption is essential. 

\begin{proof}The sky-scraper sheaf $\cE=\cO_Y$ along the zero section defines an object in $D^b_T(X,W)$; it's equivalent to the matrix factorization $\pi^*\cL^\vee \stackrel{p}{\to} \cO \stackrel{f}{\to} \pi^*\cL^\vee(\chi)$. The endomorphism algebra of $\cE$ relative to $Y$ is $\cO_Z$ so it gives a functor
$$\hom(\cE, -) : D^b_T(X, W) \to D^b_T(Z)$$ 
(see Remark \ref{rem.MFs}(6)). In general this functor does not have an adjoint. However, if we replace $X$ with its completion $\widehat{X}$ along $Y$ then it gives us an equivalence between $D^b_T(Z)$ and $D^b_T(\widehat{X}, W)$. 

The category of matrix factorizations is supported on $\Crit(W)$ so our assumption guarantees that the restriction $D^b_T(X, W)\to D^b_T(\widehat{X}, W)$ is an equivalence. The result follows.
\end{proof}

\begin{rem}\hspace{1pt}\label{rem.KP}
 \begin{enumerate}\setlength{\itemsep}{7pt}

\item A typical application of this result is when $T=\C^*$ acting trivially on $Y$, with $\chi=1$. In this case $D^b_\alpha(Z)$ is the usual derived category, and the assumption in the theorem is obviously satisfied.

\item The result generalizes to vector bundles of higher rank immediately, but we'll only need the line-bundle case. It also generalizes immediately to the case where $Y$ is an orbifold. 

\item We can formulate a $\Z_2$-graded version of this statement without choosing a $T$-action or an R-charge, by replacing $D^b(Z)$ with the $\Z_2$-graded derived category. This version is false if $Z$ is singular (see Remark \ref{rem.KPfails}).

\item We can generalize the result by allowing a non-zero superpotential on $Z$, \emph{i.e.} by setting $W=fp + g$ for some $g\in \Gamma(Y,  \cO(\chi))$ and replacing $D^b_\alpha(Z)$ with $D^b_\alpha(Z, g)$. Then the object $\cE$ lies in $D^b_T(X, fp)$ so $\hom(\cE, - )$ does indeed produce objects in $D^b_T(Z, g)$. 

\end{enumerate}
\end{rem}

The second key property of matrix factorizations is their behaviour under variation of GIT. 

\begin{thm}\cite{segal, BFK} \label{thm.VGIT} Let $V=\C^n$ be a vector space and $T'\to SL(V)$ a torus action with trivial determinant. Let $X_+$ and $X_-$ be any two generic GIT quotients $V\sslash T'$. Choose a subtorus
$$T \subset (\C^*)^n/T'$$
of the natural torus acting on $X_\pm$, and a $T'$-invariant function $W$ on $V$ transforming in a character $\chi$ of $T$. Then we have an equivalence
$$\Phi: D^b_T(X_+, W) \isoto D^b_T(X_-, W) $$
and hence an equivalence $D^b_\alpha(X_+, W)\isoto D^b_\alpha(X_-, W)$ for any R-charge $\alpha$. 
\end{thm}

\begin{rem}\hspace{1pt}\label{rem.VGIT}
 \begin{enumerate}\setlength{\itemsep}{7pt}

\item $X_+$ and $X_-$ are birational and $\Phi$ is the identity away from the exceptional loci.

\item An important special case for us is when $T'=\C^*$ is rank 1, so $X_\pm$ are the only two generic quotients. In this case there are a family of possible equivalences $\Phi$ indexed by $\Z$, and related to each other by the Picard groups of $X_\pm$.

 In general there is a large set of these equivalences $\Phi$. 

\item If we remove the $T$-actions and the R-charge then we have the same equivalence of $\Z_2$-graded categories.

\end{enumerate}
\end{rem}

\subsection{The mirror to the pair-of-pants}\label{sec.POPmirror}

Let $\Sigma$ be the pair-of-pants, a genus-zero surface with three punctures. The mirror to $\Sigma$ is the category of matrix factorizations of the superpotential $xyz$, there is an equivalence of $\Z_2$-graded categories
$$\Fuk(\Sigma) \cong D^b(\A^3, xyz) $$
\cite{AAEKO}. The category $D^b(\A^3, xyz)$ has three obvious objects which are the sky-scraper sheaves along the three co-ordinate hyperplanes $\cO/x, \cO/y$ and $\cO/z$.  Under mirror symmetry these map to arcs joining each pair of punctures. Note that there is a slight asymmetry here; on the B-side these objects are canonical, but on the A-side there is no canonical way to choose three such arcs even up to homotopy.
\pgap

To get $\Z$-graded categories we must pick a line field $\eta$ on $\Sigma$. This line field is characterized by its three winding numbers $\alpha_x, \alpha_y, \alpha_z$ around the three punctures, and they satisfy $\alpha_x + \alpha_y + \alpha_z = 2$ (see Section \ref{sec.linefields}). On the mirror side these three numbers provide a grading of the ring $\C[x,y,z]$ with $|xyz|=2$. 

We define the derived category of this graded ring following Section \ref{sec.Rcharge}. We set $T=(\C^*)^3$ acting on $\A^3$ in the obvious way, so $W=xyz$ is $\chi$-semi-invariant for $\chi=(1,1,1)$. Then $\alpha=(\alpha_x, \alpha_y, \alpha_z)$ is an R-charge and we have a $\Z$-graded category $D^b_\alpha(\A^3, xyz)$. 
\pgap

We can get a second description of this category using Kn\"orrer periodicity, Theorem \ref{thm.KP}. In the notation of that theorem we set  $Y=\A^2$ and let $T\actson Y$ with the third factor acting trivially. Then $X=\A^3$ is the total space of the line-bundle $\cL = \cO(0,0,1)$. Also $f=xy$ is a section of $\cL^\vee(\chi) $ and $W=fz$. The zero locus of $f$ is the node
$$Z = \{xy=0\} \subset \A^2$$
and we get an equivalence:
\beq{eq.KPT}\Hom(\cO/z, - ): D^b_T(\A^3, xyz) \isoto D^b_T(Z)\eeq
Note that target here is equivalent to the equivariant derived category for the rank 2 torus action $(\C^*)^2\actson Z$ (Remark \ref{rem.MFs}(4)). Then after collapsing the grading with $\alpha$ we get:
\beq{eq.KP} D^b_\alpha(\A^3, xyz) \isoto D^b_\alpha(Z)  \eeq
This equivalence maps the object $\cO/z$ to the structure sheaf $\cO_Z$. It maps the objects $\cO/x$ and $\cO/y$ to the sky-scraper sheaves $\cO_Z/x, \cO_Z/y$ along the two corresponding branches of the node. 
\pgap

There's a localization functor 
$$D^b_\alpha(\A^3, xyz) \To D^b_\alpha(\A^3_{x\neq 0},  xyz)$$
given by restricting to the open set $\{x\neq0\}$. In this open set the critical locus of the superpotential is smooth; it's the punctured affine line $\{y=z=0\}$. So Theorem \ref{thm.KP} gives us an equivalence:
\beq{eq.KPcylinder} D^b_\alpha(\A^3_{x\neq 0},  xyz) \isoto D^b_\alpha(\C[x^{\pm 1}])\eeq
Combining this with restriction gives us a functor:
\beq{eq.KPlocalized}D^b_\alpha(\A^3, xyz) \To D^b_\alpha(\C[x^{\pm 1}])\eeq
Obviously we have similar functors for the variables $y$ and $z$.  On the mirror side these functors corresponds to the inclusion of one of the cylindrical ends of $\Sigma$.

\begin{rem}
An essential observation for this paper is that these functors are \emph{not} canonical. The ambiguity is in the equivalence \eqref{eq.KPcylinder} because we can either use the functor $\Hom(\cO/z, -)$ or the functor $\Hom(\cO/y, -)$; either $y$ or $z$ can play the role of the variable $p$ in Theorem \ref{thm.KP}. On the open set $\{x\neq 0\}$ we have an isomorphism
$$\cO/y \; \cong \; (\cO/z)[\alpha_z  -1] $$
so our two equivalences differ by a shift.

In fact for our purposes even these two choices are not privileged and we will regard all shifts of these functors as equally valid. 
\end{rem}

\begin{rem}\label{rem.KPambiguous}
Even if we work with $\Z_2$-graded categories this ambiguity in Kn\"orrer periodicity is present since $\cO/y$ and $\cO/z$ differ by a shift.\footnote{Remark 3.5 in \cite{PS} is mistaken on this point.} But we will see a way around this in Section \ref{sec.Z2gradedcase}. 
\end{rem}


\begin{rem}\label{rem.localizationfornodes} If we replace $D^b_\alpha(\A^3, xyz)$ with the node $D^b_\alpha(\C[x,y]/xy)$ then the functor \eqref{eq.KPlocalized} is just restriction to the open set $\{x\neq 0\}$, which appears to be canonical. But there is still an ambiguity because we could have passed through the node $D^b_\alpha(\C[x, z]/xz)$ instead.

The node description makes the localization functors to $D^b_\alpha(\C[x^{\pm 1}])$ and to $D^b_\alpha(\C[y^{\pm 1}])$ more obvious but it obscures the third functor. To understand it, observe that under \eqref{eq.KP} the object $\cO/z$ maps to the structure sheaf on the node. Localizing to $\{z\neq 0\}$ is exactly localizing at the subcategory generated by $\cO/z$, so after Kn\"orrer periodicity it corresponds to localizing at the category of perfect complexes. So the third localization functor is:
$$D^b_\alpha(\C[x,y]/xy) \To D_{sg, \alpha}(\C[x,y]/xy) \cong D^b_\alpha(\C[z^{\pm 1}])$$
In the case $\alpha_x=\alpha_y=0$ the latter equivalence is very well-known, it says this singularity category is just the 2-periodic category of vector spaces (see Remark \ref{rem.punctureA1}). 
\end{rem}

We end this section by setting up some notation for the next one. Fix two elements $\beta_{xy}, \beta_{xz}\in \Z_{\alpha_x}$, such that $\beta_{xz} = \beta_{xy} + \alpha_z -1$. Note the symmetry here;  this is equivalent to $\beta_{xy} = \beta_{xz} + \alpha_y - 1$ because:
$$\alpha_y-1 \equiv 1-\alpha_z \mod \alpha_x$$ 
We will write $K(\beta_{xy}, \beta_{xz})$ for the equivalence
\beq{eq.betaequivalence}K(\beta_{xy}, \beta_{xz})\!:\; D^b_\alpha(\A^3_{x\neq 0},  xyz) \isoto D^b_\alpha(\C[x^{\pm 1}])\eeq
given by taking Hom's from the object:
$$(\cO/y)[-\beta_{xz}] \cong (\cO/z)[-\beta_{xy}]$$
This equivalence sends:
$$\cO/y \mapsto \cO[\beta_{xz}] \aand \cO/z \mapsto \cO[\beta_{xy}] $$
\pgap

\begin{rem} \label{rem.KPfails} The example of this section sheds some light on both (a) the technical assumption in Theorem \ref{thm.KP} and (b) our rather convoluted definition of $D^b_\alpha(Z)$. 

Recall that Kn\"orrer periodicity maps the object $\cO/x \in D^b_\alpha(\A^3, xyz)$ to the sky-scraper sheaf $\cO_Z/x$ along one branch of the node. The object $\cO/x$ has endomorphism algebra:
$$\C[y,z]/yz$$
This is true for any $\alpha$, and also in the $\Z_2$-graded setting. But the endomorphisms of the sheaf $\cO_Z/x$ on $Z$ are a little more subtle. If we use the naive derived category of the graded ring $\C[x,y]/xy$ - the category we called $D^b_\alpha(Z)^{obv}$ in Section \ref{sec.Rcharge} - then we get one of two possible answers:
\al{(i) && \C[y,z]/yz  \hspace{1cm}\mbox{if } \alpha_z\neq 0  \\
(ii) && \C[y][[z]]/yz \hspace{1cm}  \mbox{if } \alpha_z =  0  }
In some sense $(ii)$ is the `correct' answer; in case $(i)$ we only get polynomials because - by definition - the Hom spaces in a graded category are the direct sum of their graded pieces. We also get answer $(ii)$ if we do the calculation in the $\Z_2$-graded derived category. This demonstrates that Kn\"orrer periodicity would have failed if we'd used either $D^b_\alpha(Z)^{obv}$ and set $\alpha_z=0$, or if we'd used the $\Z_2$-graded derived category.\footnote{The version in \cite[Prop. 2.3]{Nadler} uses not the usual $\Z_2$-graded derived category but rather the $\Z$-graded derived category with the grading collapsed. This is similar to our $D^b_\alpha(Z)$.} It also demonstrates that our $D^b_\alpha(Z)$ is different from $D^b_\alpha(Z)^{obv}$ in the case $\alpha_z=0$. 

Another closely-related way that we can break Kn\"orrer periodicity is to use a smaller torus action. Suppose instead of our rank 3 torus $T$ we just use the rank 2 torus $T'=(\C^*)^2$ acting on the $x,y$ co-ordinates. Then the assumption in Theorem \ref{thm.KP} fails; $\Crit(xyz)$ includes the $z$-axis, and $T'$ acts trivially on this component. The problem is that sky-scraper sheaves at points up the $z$-axis define non-zero objects of $D^b_{T'}(\A^3, xyz)$ which are orthogonal to $\cO/z$, and are hence in the kernel of \eqref{eq.KPT}.

Finally this gives another demonstration of why we need to take pre-triangulated closure, similar to Example \ref{eg.localisingA1}. If $\alpha_z=0$ then the category $D^b_\alpha(Z)$ contains cones of the form 
$$\cO_Z/x \stackrel{z-\lambda}{\To} \cO_Z/x$$
 which do not exist in $D^b_T(Z)$. In $D^b_\alpha(\A^3, xyz)$ these correspond to point sheaves up the $z$-axis as discussed just above. Rather curiously these objects in $D^b_\alpha(Z)$ are orthogonal to the structure sheaf, so they are acyclic but non-zero. 
\end{rem}

\section{Categories from decorated trivalent graphs}\label{sec.categoriesfromgraphs}

In this section we set up some combinatorial abstractions based on our previous discussions of the pair-of-pants and it's mirror.

\begin{defn} Let $E$ be a set with three elements. A \emph{decoration} of $E$ is the following data:
\begin{enumerate}\item A function $\alpha: E \to \Z$,  $ x \mapsto \alpha_x$, satisfying:
$$\sum_{x\in E} \alpha_x = 2 $$
\item For each distinct ordered pair $(x,y)\in E\times E$ an element $\beta_{xy} \in \Z_{\alpha_x}$, satisfying
$$\beta_{xz} = \beta_{xy} + \alpha_z -1 $$
for each distinct triple $(x,y,z)$. 
\end{enumerate}
\end{defn}
If $E=\{x,y,z\}$ it's enough to specify $\alpha_x$ and $\alpha_y$ since then $\alpha_z$ is determined; similarly it's enough to specify (say) $\beta_{xy}, \beta_{yz}$ and $\beta_{zx}$ since then the other three $\beta$ elements are determined. But we can't do this without breaking the symmetry somehow. 

 Given a decoration of $E$ we can define
\beq{eq.gammas}\gamma_{xy}= \beta_{xy} - \beta_{yx} \;\in \Z/(\alpha_x, \alpha_y)\eeq
for each distinct ordered pair $(x,y)$. Obviously these invariants are anti-symmetric, they also obey the relation
$$\gamma_{xy} + \gamma_{yz} + \gamma_{zx}  = 1 \; \in \Z/(\alpha_x, \alpha_y, \alpha_z) $$
(note that this group is either trivial or $\Z_2$).  If we add 1 to each $\beta_{xy}$ then we get a new decoration and obviously the $\gamma$'s do not change. Conversely if two decorations have the same $\alpha$ and the same $\gamma$ then it's elementary to prove that their $\beta$'s are also the same up to the addition of a constant.
\pgap 

 By a \emph{graph} we mean the usual thing, with the conventions that we allow loops, multiple edges between vertices, and external edges. We will emphasise the role of half-edges, so an internal edge is the union $x\cup y$ of two half-edges $x$ and $y$, and an external edge is a single half-edge. We'll sometimes write a vertex as the union $x\cup y \cup z$ of the half-edges that meet it. For a graph $\Gamma$ we'll write $\partial \Gamma$ for the set of external edges. 
 
 We'll be interested almost excusively in trivalent graphs with $\partial \Gamma \neq \emptyset$.

\begin{defn} Let $\Gamma$ be a trivalent graph. A \emph{decoration} of $\Gamma$ is a choice, for each vertex, of a decoration $(\alpha, \beta)$ on the set of half-edges at that vertex. For each internal edge $x\cup y$ we require that:
$$\alpha_{x}  + \alpha_{y} = 0 $$
\end{defn}

So our decoration assigns an integer $\alpha_x$ to every half-edge in $\Gamma$ and an element $\beta_{xy}\in \Z_{\alpha_x}$ to every ordered pair of half-edges which meet the same vertex. We don't require any relationship between the $\beta$'s at different vertices.

Summing $\alpha$ over all half-edges gives
\beq{eq.eulerchar} \sum_{x\in \partial \Gamma} \alpha_x = 2v \eeq
where $v$ is the number of vertices. So it's impossible to decorate a graph with no external edges; conversely it's trivial to show that a graph with external edges does admit a decoration.  Note that if $\Gamma$ is connected then $2v= 4g +2e - 4$ where $g$ is the genus and $e$ is the number of external edges. 

 Another elementary observation is that the set of possible $\alpha$'s (if $\partial\Gamma\neq \emptyset$) forms a torsor for the relative homology group $H_1(\Gamma, \partial\Gamma)$. If $\Gamma$ is connected this is $\Z^{g+e-1}$. If we fix the values for $\alpha$ on the external edges subject to \eqref{eq.eulerchar} then the remaining possibilities form a torsor for $H_1(\Gamma)$. 
\pgap

Given a decorated trivalent graph we are going to construct an associated category.  We will need the following data:
\begin{enumerate} \item For each vertex $v$ a dg-category $\cC_v$.
\item For each half-edge $x$ a dg-category $\cC_x$.
\item For each internal edge $x\cup y$ an equivalence $\cC_x \cong \cC_y$.
\item For each half-edge $x$, meeting a vertex $v$, a functor from $\cC_v$ to $\cC_x$. 
\end{enumerate}
This data gives us a diagram of dg-categories, and then we can take the homotopy limit of the diagram. Note that since $\Gamma$ is just a 1-dimensional simplicial set we don't need any higher-categorical data to take this limit. 

For a decorated trivalent graph $(\Gamma, \alpha, \beta)$ we choose the data as follows:
\begin{enumerate}
\item Let $v=x\cup y\cup z$ be a vertex. We set $\cC_v$ to be $D^b_\alpha(\A^3, xyz)$.
\item For a half-edge $x$ we set $\cC_x$ to be $D^b_\alpha(\C[x^{\pm 1}])$.
\item For an internal edge $x \cup y$ we use the equivalence $\cC_x\cong \cC_y$ induced by the isomorphism of graded rings $\C[x^{\pm 1}] \isoto \C[y^{\pm 1}]$ which maps $x\mapsto y^{-1}$.
\item For a half-edge $x$ meeting the vertex $v=x\cup y\cup z$, we have  functor $\cC_v$ to $\cC_x$ given by localizing $\cC_v$ at the variable $x$, then applying the equivalence $K(\beta_{xy}, \beta_{xz})$ from \eqref{eq.betaequivalence}.
\end{enumerate}
We'll denote the limit category by 
$$\cC_\Gamma$$
though it does depend on the decoration $(\alpha, \beta)$ as well as the graph $\Gamma$. 
\pgap

Notice that at each external edge $x$ we have a functor from $\cC_\Gamma$ to $\cC_x$, this is extra data beyond the category $\cC_\Gamma$ itself, which we can choose either to keep or forget. 

In this paper we will often have two decorated graphs $\Gamma$ and $\Gamma'$ with the same sets of external edges, \emph{i.e.} with a given bijection $\partial \Gamma \isoto \partial \Gamma'$, and with $\alpha_x = \alpha'_{x}$ for each $x\in \partial \Gamma$.  In this situation we will say that the two categories $\cC_\Gamma$ and $\cC_{\Gamma'}$ are \emph{equivalent relative to their external edges} if we have an equivalence
$$\cC_\Gamma \isoto \cC_{\Gamma'} $$
which intertwines the functors 
$$\cC_\Gamma \to \cC_x \aand \cC_{\Gamma'} \to \cC_x $$
for each external edge $x$. 

\begin{rem}\label{rem.externalring}
A related point is the observation that the homotopy category of $\cC_\Gamma$ is linear over the ring $R = \C[x_1,..., x_k]/(x_ix_j)$ where $x_i$ are the external edges (though it's not immediately obvious if this lifts to the dg-level). But saying that an equivalence $\cC_\Gamma \isoto \cC_{\Gamma'}$ is R-linear is not quite as strong as saying it's an equivalence relative to external edges; the difference is shift functors at the external edges.
\end{rem}

\subsection{Trivial modifications}\label{sect.trivialmods}

The $\beta$ part of the decoration has some obvious redundancies, there are several ways we can change it without affecting the category $\cC_\Gamma$. 

\begin{defn}\label{defn.trivmod} A \emph{trivial modification} of a decorated graph is one of the following three operations:
\begin{enumerate}\setlength{\itemsep}{5pt}

\item[(V)] At a vertex $x\cup y \cup z$ we can add an integer $n$ to all six values $\beta_{xy}, \beta_{xz},..., \beta_{zy}$. This corresponds to acting with the autoequivalence $[-n]$ on $\cC_v$.

\item[(I)] Suppose we have an internal edge $x_1\cup y_1$ joining vertices with half edges $\{x_1, x_2, x_3\}$ and $\{y_1, y_2, y_3\}$. We can add an integer $m$ to $\beta_{x_1x_2}$, $\beta_{x_1x_3}$,  $\beta_{y_1y_2}$ and $\beta_{y_1y_3}$. This corresponds to acting with the autoequivalence $[m]$ on the categories $\cC_{x_1}$ and $\cC_{y_1}$ simultaneously.
\item[(E)] Suppose we have a vertex $x\cup y\cup z$ and $x$ is an external edge. We can add an integer $m$ to $\beta_{xy}$ and $\beta_{xz}$. This corresponds to acting with the autoequivalence $[m]$ on $\cC_{x}$. 
\end{enumerate}
\end{defn}

So the category $\cC_\Gamma$ only depends on the equivalence class of $\beta$ under trivial modifications.

Modifications of types (V) and (I)  preserve the category $\cC_\Gamma$ `relative to its external edges', \emph{i.e.} they commute with the functors $\cC_\Gamma\to \cC_x$ at each external edge $x$. Modifications of type (E)  do not have this property. We will sometimes refer to the first two types as \emph{relative trivial modifications}.

\begin{rem}\label{rem.gammadecoration}
 If we take the orbit of $\beta$ under trivial modification of type (V) then this is the same as remembering only the $\gamma$ invariants \eqref{eq.gammas} at each vertex; so $\cC_\Gamma$ only depends on these $\gamma$ invariants. 

It follows that if we decorate a trivalent graph by $\alpha$ data and $\gamma$ data (instead of $\beta$ data) then the category $\cC_\Gamma$ is well-defined. Note that this kind of decorated graph is exactly what arose in our discussion of line fields on pair-of-pants decompositions in Section \ref{sec.linefields}.
\end{rem}

\begin{rem} \label{rem.deltadecoration}Now let's consider the orbit of $\beta$ under trivial modifications of type (I). Say we have an internal edge  $x_1\cup y_1$, joining vertices with half edges $\{x_1, x_2, x_3\}$ and $\{y_1, y_2, y_3\}$ (the two end vertices need not be distinct). We can define four invariants
$$\delta_{x_iy_j} = \beta_{x_1 x_i} - \beta_{y_1 y_j} \; \in \Z/\alpha_{x_1}\;\; i,j\in[2,3] $$
(which makes sense since $\Z/\alpha_{y_1}=\Z/\alpha_{x_1}$). They satisfy relations of the form
$$\delta_{x_2y_3} = \delta_{x_2y_2} - \alpha_{y_3} + 1$$
so any one of them determines the other three, and it's immediate that replacing $\beta$'s by $\delta$'s exactly classifies the orbits under trivial modifications of type (I) and (E). So we could also decorate our graphs by $\alpha$ and $\delta$ data and still define the cateory $\cC_\Gamma$. 

Also note that a trivial modification of type (I) corresponds on the A-side to moving the marked point on that cuff, see Remark \ref{rem.movep}.  
\end{rem}

\begin{lem}\label{lem.genus0trivmod} Let $\Gamma$ have genus zero. Then for a given $\alpha$ any two decorations $(\alpha, \beta)$ and $(\alpha, \beta')$ are equivalent under trivial modifications. 
\end{lem} 
\begin{proof} On each connected component pick a vertex and then walk away from it along internal edges, applying operations (I) then (V) on each internal edge. Finally apply (E) on the external edges.   
\end{proof}
So in genus zero the $\beta$ data is actually irrelevant for the category $\cC_\Gamma$. In higher genus this is not true because summing the values of $\beta$ around loops in the graph gives us invariants.

Suppose we have an oriented cycle $c\subset \Gamma$ consisting of edges $x_0\cup y_1$, $x_1\cup y_2$, ..., $x_k \cup y_0$, where for each $i$ the half-edges $x_i$ and $y_i$ meet at a vertex.  Let's define $I_c$ to be the ideal:
$$I_c = (\alpha_{x_0}, \alpha_{y_1}, ..., \alpha_{x_k}, \alpha_{y_0})  = (\alpha_{x_0}, ..., \alpha_{x_k}) \; \subset \Z $$
Then the quantity

\begin{align}\label{eq.bc} b_c &= \beta_{x_0y_0} -\beta_{y_0x_0} + \beta_{x_1y_1} - ... -\beta_{y_0x_k}\quad \quad \in \Z/I_c \nonumber\\
& = \gamma_{x_0y_0} + .... + \gamma_{x_k y_k}\nonumber \\
& = \delta_{y_0 x_1} + ... + \delta_{y_k x_0}  \end{align}
is invariant under trivial modifications (the second expression is a sum over vertices, the third is a sum over edges).

\subsection{The $\Z_2$-graded case}\label{sec.Z2gradedcase}

The above construction is designed to mirror Fukaya categories of surfaces with line fields, so they're $\Z$-graded. If we only want our Fukaya categories to be $\Z_2$-graded then one can use a `mod 2 line field' which means a certain kind of element of $H^1(\P T\Sigma, \Z_2)$; see \cite{LPinvs}. In our mirror construction this simply corresponds to taking each $\alpha_i$ to be an element of $\Z_2$ instead of an integer, and each $\beta_{ij}$ to be an element of $\Z/(2, \alpha_i, \alpha_j)$. 

However, unlike the $\Z$-graded case, there is a canonical choice of `mod 2 line field' induced from the orientation of $\Sigma$. This is because there is a canonical isomorphism from the set of `mod 2 line fields' to $H^1(\Sigma, \Z_2)$ given by taking winding numbers in $\Z_2$. For the the canonical choice all winding numbers are zero and all Lagrangians are gradeable; this is why one can define $\Z_2$-graded Fukaya categories without mentioning line fields.

On our mirror side we can make our $\alpha$ values lie in $\Z_2$ instead of $\Z$, then it's obvious that the set of possible $\alpha$'s is canonically the relative homology group $H_1(\Gamma, \partial\Gamma; \Z_2)$. In particular the canonical choice is just to set $\alpha_x=0\in \Z_2$ for all $x$, then at each vertex we have the usual $\Z_2$-graded category $D^b(\A^3, xyz)$. But we still need to choose the $\beta_{ij}$'s, \emph{i.e.} we need to choose the following data:

\begin{defn} A \textit{weak decoration} of a trivalent graph $\Gamma$ is a choice of element $\beta_{xy}\in \Z_2$ for each ordered pair of half-edges that meet at a vertex, such that at any vertex $x\cup y\cup z$ we have $\beta_{xz}= \beta_{xy} + 1$. 
\end{defn}

As before we only need the equivalence class of a weak decoration up to trivial modifications. It's easy to see that the set of weak decorations up to trivial modifications forms a torsor for $H^1(\Gamma, \Z_2)$. It's slightly less obvious that this torsor is canonically trivial.
\begin{lem}\label{lem.H1mod2} Let $\Gamma$ be a trivalent graph. The set of weak decorations of $\Gamma$, taken up to trivial modifications, is canonically isomorphic to $H^1(\Gamma, \Z_2)$.
\end{lem}
\begin{proof} For a given weak decoration $\beta$ we define a function on cycles in $\Gamma$, by sending 
$$c \mapsto b_c \; \in \Z_2$$
as in \eqref{eq.bc}. We claim that $b_c$ is actually a cocycle, \emph{i.e.} if $c=c_1 + c_2$ then $b_c = b_{c_1}+ b_{c_2}$. These two expressions are not identical by definition, they differ by 1 for every vertex where $c_1$ and $c_2$ meet or diverge. But there is an even number of such vertices.
\end{proof}

\begin{rem} This lemma is a combinatorial counterpart to the result for surfaces mentioned above. If a surface $\Sigma$ has a pair-of-pants decomposition encoded by $\Gamma$ then there is a short exact sequence:
\beq{SES}0\to H^1(\Gamma) \to H^1(\Sigma) \to H_1(\Gamma, \partial \Gamma) \to 0\eeq
Indeed the Mayer-Vietoris sequence for the pair-of-pants decomposition gives an exact sequence
$$ \Z^v \to \Z^i \to H^1(\Sigma) \to \Z^{2v} \to \Z^i \to 0 $$
where $v$ and $i$ are the numbers of vertices and internal edges of $\Gamma$. The cokernel of the first map is obviously $H^1(\Gamma)$. If we add a bivalent vertex to every internal edge of $\Gamma$ then we see that $H_1(\Gamma, \partial\Gamma)$ is the kernel of a map $\Z^{3v} \to \Z^{v+i}$, and this agrees with the kernel of the final map in the sequence above.

If we take the sequence \eqref{SES} with $\Z_2$ coefficients then the middle term is the set of `mod 2 line fields' on $\Sigma$, and the final term is the choice of `mod 2 $\alpha$ data'. So if $\alpha\equiv 0$ then the remaining choices lie in $H^1(\Gamma, \Z_2)$.
\end{rem}

Since a trivalent graph $\Gamma$ has a canonical weak decoration (up to trivial modifications), there is a canonical $\Z_2$-graded category associated to it. We'll denote this category by $\cD_\Gamma$ to distinguish it from the $\Z$-graded version. By \cite{Lee, PS} it is equivalent to the Fukaya category of the associated punctured surface. 

\begin{rem}\label{rem.PSplanar} This discussion fills a rather small gap in the recent paper \cite{PS2}. There the authors give a very similar construction of a $\Z_2$-graded category to a trivalent graph, but it uses the erroneous claim that Kn\"orrer periodicity is canonical (see Remark \ref{rem.KPambiguous}).  In the earlier \cite{PS} they only consider only \emph{planar} trivalent graphs and give a slightly more detailed construction which avoids this gap. On a planar graph there is a more obvious construction of the canonical weak decoration (see Remark \ref{rem.planarweakdecoration}).
\end{rem}

\section{Geometric mirrors with R-charge and twists}

In this section we discuss situations in which the category $\cC_\Gamma$ has a geometric interpretation, and also some algebro-geometric equivalences between them.

\subsection{One-dimensional mirrors}\label{sec.STZ}

\begin{figure}
\begin{tikzpicture}[scale =1]
\draw (0,0) --(4,0);
\draw (1,0)--(1, 1);
\draw (3,0)--(3,1);
\draw[dashed] (4,0)--(6,0);
\draw (6,0)--(8,0);
\draw (7,0)--(7,1);
\node [below] at (0,0) {$x_1$};
\node [above] at (1,1) {$z_1$};
\node [below] at (1.5,0) {$y_1$};
\node [below] at (2.5,0) {$x_2$};
\node [above] at (3,1) {$z_2$};
\node [below] at (3.5,0) {$y_2$};
\node [below] at (6.5,0) {$x_k$};
\node [above] at (7,1) {$z_k$};
\node [below] at (8,0) {$y_k$};
\end{tikzpicture}
\caption{A straight tree.}
\label{fig.straighttree}
\end{figure}
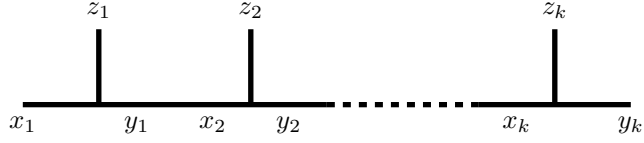

Let $\Gamma$ be a trivalent graph with the property that every vertex meets an external edge. Such a $\Gamma$ can only have genus zero or one. For genus zero it must be a `straight tree' as pictured in Figure \ref{fig.straighttree}. For genus one it must be a `wheel', obtained from a straight tree by gluing $x_1$ to $y_k$. So our surface $\Sigma$ is either a punctured sphere or a punctured torus. 

Equip $\Gamma$ with a decoration $(\alpha, \beta)$. In this situation the category $\cC_\Gamma$ has a geometric interpretation as the (perhaps twisted) derived category of a nodal curve $C$. The untwisted version of this is well-known, e.g.~\cite{STZ}. 
\pgap

At each vertex we have the category $D^b_\alpha(\A^3, x_iy_iz_i)$, but we can apply Theorem \ref{thm.KP} and replace this with the derived category of a node $Z_i = \{x_iy_i=0\}$. Gluing these nodes together produces a chain $C$ of rational curves, meeting at nodes. In the genus one case these rational curves form a loop, and in genus zero they form a line with an $\A^1$ at each end. The former is an example of a `balloon ring' and the latter a `balloon chain' latter is an example of a `balloon chain'  \cite{STZ}.

We have a torus action $T'=(\C^*)^r\actson C$ where $r$ equals the number of irreducible components, it acts by rotating each $\P^1$ and $\A^1$.  Our $\alpha$ data specifies a 1-parameter subgroup $\C^* \to T'$. This defines our category $D^b_\alpha(C)$, by using $\alpha$ to collapse the grading on the $T'$-equivariant derived category (Remark \ref{rem.Dbalpha}(2)). 

Since equivariant derived categories form a sheaf (\emph{e.g.}~\cite{Toen}) we can reconstruct $D^b_\alpha(C)$ as the limit of a diagram in the same way that we constructed $\cC_\Gamma$. To spell this out: the open subsets $\{x_i \neq 0 \}$ and $\{y_i\neq 0\}\subset Z_i$ are copies of $\A^1\setminus 0$ and we have restriction functors from $D^b_\alpha(Z_i)$ to each of them. We have also have equivalences $D^b_\alpha(\{y_i\neq 0\}) \isoto D^b_\alpha(\{x_{i+1}\neq 0\})$ since we've glued the intersections together. This gives a diagram and $D^b_\alpha(C)$ is the limit.

Comparing this with the definition of $\cC_\Gamma$ we see that 
 $$D^b_\alpha(C)\cong \cC_\Gamma$$
for a \emph{specific choice} of $\beta$ data. At vertex $i$ the object $\cO/z_i\in D^b_\alpha(\A^3, x_iy_iz_i)$ maps to the structure sheaf on $Z_i$ which then restricts to the structure sheaf on each intersection. For a general $\beta$ this object would map to some shift of the structure sheaf on the intersections, see \eqref{eq.betaequivalence}. So to get $D^b_\alpha(C)$ we must set
$$\beta_{x_iy_i} = 0 \aand \beta_{y_i x_i} = 0 $$
for all $i$ (and hence $\beta_{x_iz_i} =\alpha_{z_i}-1 $ and $\beta_{y_iz_i} = \alpha_{z_i}-1$).
\pgap 

In genus zero the $\beta$ data does not matter (Lemma \ref{lem.genus0trivmod}) so in fact $\cC_\Gamma \cong D^b_\alpha(C)$ for any choice of $\beta$. In particular if we set $\alpha_{x_i}=\alpha_{y_i}=0$ (so $\alpha_{z_i}=2$) for each $i$ then $\cC_\Gamma$ is the usual derived category $D^b(C)$. 
\pgap

The genus one case is more interesting. Recall from Section \ref{sect.trivialmods} that we can extract an invariant from the $\beta$ data of our decoration, it's the sum:
\al{ b &= \beta_{x_1y_1}-\beta_{y_1x_1} + ... + \beta_{x_ky_k}-\beta_{y_kx_k} \\
& = \delta_{y_2 x_1} + ... + \delta_{y_k x_{k-1}} + \delta_{y_1 x_k} \;\in \Z/(\alpha_{x_1}, ..., \alpha_{x_k}) }
This residue is invariant under trivial modifications; in fact it's easy to prove that for a given $\alpha$ it classifies the possible $\beta$'s up to trivial modifications. The choice of $\beta$ that produces $D^b_\alpha(C)$ has $b=0$. 

What happens if we choose a different $\beta$? By construction our curve $C$ comes with an open cover $\{U_1,..., U_k)\}$ in which each triple overlap is empty. Each term $\delta_{y_j x_{j-1}}$ in $b$ corresponds to the overlap $U_j\cap U_{j-1}$ and if we lift each to an integer $\tilde{\delta}_{y_j x_{j-1}}$ then we have a cocycle, hence a class
$$\tilde{b}  \in  H^1_{\acute{e}t}(C, \Z)= \Z$$
(with $[\tilde{b}]=b$). This class provides a twisting of the category $D^b_\alpha(C)$ as described in Section \ref{sec.H1twists} and  it's evident that
$$ \cC_{\Gamma} \cong D^b_\alpha(C, \tilde{b}) $$
is the twisted derived category.
\pgap

 Note that $\alpha$ and $\tilde{b}$ are completely independent here; depending on the decoration chosen we can get the untwisted category $D^b_\alpha(C)$ with R-charge, or the twisted version $D^b(C, \tilde{b})$ of the usual derived category without R-charge. In general we get both R-charge and a twisting. 

\begin{rem}\label{rem.bmodalpha}
It's also apparent that we have an equivalence
$$D^b_\alpha(C, \tilde{b}) \isoto D^b_\alpha(C, \tilde{b} + \alpha_{x_i}) $$
for any $i$, since the category $\cC_\Gamma$ only depends on the residue $b$.  This can be understood as tensoring by a `line-bundle'. Consider the line-bundle which has degree 1 on the $i$th component of $C$ and degree zero everywhere else, so the transition function is $x_i$ on $U_i\cap U_{i+1}$ and the identity elsewhere. Since $x_i$ is a variable of degree $\alpha_{x_i}$ these transition functions define an object in $D^b_{\alpha}(C, \alpha_{x_i})$, and tensoring by this object gives the equivalence above. 

On the mirror side this is a Dehn twist, it provides an equivalence between the Fukaya categories of $\Sigma$ for two different line fields.
\end{rem}

The case $k=1$ is particulary interesting and we'll return to it in Section \ref{sec.nodalelliptic}.

\subsection{Hori-Vafa mirrors}\label{sec.HoriVafa}

Let $X$ be a smooth toric Calabi-Yau three-fold. There's a canonical superpotential $W: X\to \C$ whose zero locus is the union of all the toric divisors; in each toric chart $\A^3\subset X$ it restricts to the function $xyz$. 

The polytope for $X$ is the cone on a polygon, which is equipped with a complete triangulation. The dual graph $\Gamma$ is trivalent and has no loops or multiple edges. One interpretation of $\Gamma$ is that it describes the critical locus of $W$, since $\operatorname{Crit}(W)$ is a curve in $X$ consisting of one rational curve for each internal edge of $\Gamma$ and one affine line for each external edge. 

It was observed in \cite{Lee} and \cite{PS} that the $\Z_2$-graded category $D^b(X,W)$ can be constructed from this graph, \emph{i.e.}
$$ \cD_\Gamma \cong D^b(X, W)$$
This follows immediately from the fact that $D^b(X,W)$ forms a sheaf on $X$. The pair $(X,W)$ is the Hori-Vafa mirror to the surface encoded by $\Gamma$ and this observation forms part of Lee and Pascaleff-Sibilla's proof of homological mirror symmetry. 
\pgap

We'll now discuss how to improve this to a $\Z$-graded equivalence. We have a torus action $T\actson X$ and $W$ is semi-invariant, so just need to choose an R-charge $\alpha: \C^* \to T$ for which $W$ has weight 2. In every toric chart this gives us a grading on the variables satisying $|xyz|=2$ so we get $\alpha$ data on the graph $\Gamma$. We claim this $\alpha$ data is of a special form. 

Say $u\cup v$ is an internal edge joining vertices $x\cup y\cup u$ and $v\cup z\cup w$, as in Figure \ref{fig.IHmoves}(a).  The planar structure gives us a preferred bijection between $\{x, y\}$ and $\{z,w\}$ according to which side of $u\cup v$  they lie; let's suppose $x\mapsto z$ and $y\mapsto w$. This data is describing two toric charts in $X$ related by a transition function
\beq{transitionfn} v = u^{-1}, \quad \quad z = x u^\sigma, \quad\quad w = y u^\tau \eeq
for some $\sigma, \tau\in \Z$ with $\sigma+\tau=2$. It follows that the $\alpha$ data must satisfy the condition:
\beq{3foldcondition}\alpha_x \equiv \alpha_z \aand \alpha_y \equiv \alpha_w \quad \mbox{mod } \alpha_u \eeq

\begin{eg}\label{eg.3foldflop} Let $X$ be the resolved conifold, the total space of $\cO(-1)^{\oplus 2}$ over $\P^1$.  We can express $X$ as a GIT quotient of $\C^4_{s, t, p, q}$ by $\C^*$ acting with weights $(1,1,-1,-1)$, where $s=t=0$ is unstable. The superpotential is $W=stpq$. Since $X$ is covered by two toric charts the associated trivalent graph is  genus zero  with two vertices. If we label this graph as in Figure \ref{fig.IHmoves}(a) then the four external edges of the graph correspond to the functions:
$$x=sp,\quad \quad  y = sq, \quad\quad z = tp \aand w=tq$$ 
 To put an R-charge on $X$ we must assign four weights $\alpha_s,\alpha_t,\alpha_p, \alpha_q$ which sum to 2. Then the external edges have weights $\alpha_x = \alpha_s + \alpha_p$ etc.~and the internal edge has weight $\alpha_u = \alpha_t-\alpha_s$. 

 For these charts the transition functions have $\sigma=\tau=1$ and correspondingly $\alpha_z = \alpha_x + \alpha_u$ and $\alpha_w = \alpha_y + \alpha_u$. This is a stronger version of condition \eqref{3foldcondition}, it's equivalent to the condition
\beq{eq.flopcondition} \alpha_x + \alpha_w = 2 \aand \alpha_y+\alpha_z=2 \eeq
and it's easy to see that any $\alpha$ data of this form can be lifted to an R-charge on $\C^4$. 

We will say more about this example in the next section.
\end{eg}

Now suppose that $\Gamma$ is any trivalent graph, embedded in a plane. Suppose we choose an $\alpha$ on $\Gamma$ satisying the condition \eqref{3foldcondition}. We claim there is now a canonical choice of $\beta$ data (up to trivial modifications).   We express this canonical choice in terms of the $\delta$ data as in Remark \ref{rem.deltadecoration}; at each edge we set: 
$$\delta_{xz}=0 \aand \delta_{yw}=0$$ Note that $\delta_{yw} = \delta_{xz} + \alpha_y-\alpha_w$, so this is only consistent because of \eqref{3foldcondition}. 

\begin{rem}\label{rem.planarweakdecoration} Recall from Lemma \ref{lem.H1mod2} that there is a canonical weak decoration on any graph (up to trivial modifications). If the graph is planar then this construction produces it.
\end{rem}
\begin{rem}\label{rem.distinguishedcycles}
 Note that since $\Gamma$ is planar there are some distinguished cycles: the boundaries of the bounded connected components of $\R^2\setminus \Gamma$. Recall that for any cycle $c$ we have our invariant $b_c$ as in \eqref{eq.bc}. If $c$ is one of the distinguised cycles and we use  our canonical choice of $\beta$ then the invariant $b_c$ will clearly be zero. 
\end{rem}

\begin{lem}\label{lem.HVmirror} Let $X$ be a smooth toric Calabi-Yau 3-fold equipped with an R-charge $\alpha$ and the canonical superpotential $W$. Let $\Gamma$ be the associated planar trivalent graph, equipped with the associated decoration $(\alpha, \beta)$. Then:
$$\cC_\Gamma \cong D^b_\alpha(X, W) $$
\end{lem}
\begin{proof} Every vertex $x\cup y \cup u \in \Gamma$ corresponds to a toric chart $X_{xyu}\subset X$ which is canonically isomorphic to $\A^3_{x,y,u}$, and the half-edge $u$ corresponds to the open subset $(X_{xyu})_{u\neq 0}$. Suppose $u$ is joined to a half-edge $v$ which meets a vertex $v\cup z\cup w$, then we have a second chart $X_{vzw}$ and the open set $(X_{vzw})_{v\neq 0}$ is the same as $(X_{xyu})_{u\neq 0}$. But the two canonical isomorphisms differ by a  transition function of the form \eqref{transitionfn}.

 Hence we can form a diagram of categories similar to the one that defines $\cC_\Gamma$. To a vertex $x\cup y\cup u$ we attach $D^b_\alpha(X_{xyu}, W)$  and to the half-edge $u$ we attach the localization of the vertex category at $u\neq 0$. If $u$ meets $v$ then the categories for $u$ and $v$ are identical. The limit over this diagram is $D^b_\alpha(X, W)$.

The category $D^b_\alpha(X_{xyu}, W)_{u^{-1}}$ attached to a half-edge $u$ is equivalent via Kn\"orrer periodicity to $D^b_\alpha(\C[u^{\pm1}])$ but this equivalence is only canonical up to shifts. If we choose one such equivalence, and $u$ meets the half-edge $v$, we get an induced equivalence $D^b_\alpha(X_{vzw}, W)_{v^{-1}}\isoto D^b_\alpha(\C[v^{\pm 1}])$. If we make one of these choices for each edge of the graph then it's clear that $D^b_\alpha(X, W)$ will be equivalent to $\cC_\Gamma$ for some choice of $\beta$; we just need to argue that the correct choice of $\beta$ is the canonical one described above.

 To see this take two adjacent vertices $x\cup y\cup u$ and $v\cup z\cup w$.  From the transition function \eqref{transitionfn} the objects $\cO/x\in D^b_\alpha(X_{xyu}, W)$ and $\cO/z \in D^b_\alpha(X_{vzw}, W)$ restrict to the same object where the charts overlap, and the same is true of $\cO/y$ and $\cO/w$. It follows that  our decoration on $\Gamma$ must satisfy $\delta_{xz}=\delta_{yw}=0$. 
\end{proof}

\begin{rem}\label{rem.twistedCrit}
Once again we should ask: what happens if we choose a different $\beta$? The distinguished cycles in our planar graph $\Gamma$ (Remark \ref{rem.distinguishedcycles}) form a basis for $H_1(\Gamma)$. For each such cycle $c$ we have our invariant $b_c$, and if we choose a lift of each $b_c$ to an integer $\tilde{b}_c$ then we get a class:
$$\tilde{b} \in H^1(\Gamma, \Z) =  H^1_{\acute{e}t}(\Crit(W), \Z) $$
Associated to this class there should be a twisted version $D^b_{\alpha}(X, W, \tilde{b})$ of the category of matrix factorizations, and it's more-or-less clear that 
$$\cC_{\Gamma} \cong D^b_{\alpha}(X, W, \tilde{b})$$
for a general $\beta$. But as far as we are aware the foundations for this twisted category have not been written down. 
\end{rem}

\subsection{Standard flops}\label{sec.flops}

Let $\Delta$ be the genus zero graph with two vertices as in Figure \ref{fig.IHmoves}. We observed in Example \ref{eg.3foldflop} that this is the graph associated to the resolved conifold $X$, so if we choose a decoration $(\alpha, \beta)$ with $\alpha$ satisfying the condition \eqref{eq.flopcondition} then we have an R-charge on $X$ and we have an equivalence:
$$\cC_\Delta \cong D^b_\alpha(X, W)$$
The choice of $\beta$ is irrelevant since $\Delta$ has genus zero. 
\pgap

If $\alpha$ does not satisfy \eqref{eq.flopcondition} then it is still possible to describe $\cC_\Delta$ as the category of matrix factorizations on a toric variety; we just have to increase the dimension from 3 to 5. Let $\C^*$ act on $\C^6$ with weights
$$\begin{matrix} s & t & p & q & y & z \\ 1 &1& -1& -1& 0& 0\end{matrix} $$
where the letters denote the co-ordinates on $\C^6$. Let  $Y$ be the GIT quotient where $s=t=0$ is unstable, so $Y$ is simply $X\times \A^2_{y,z}$. Now equip $W$ with the superpotential
$$W = tpy + sqz $$
(this is an invariant function on $\C^6$ so it descends to $Y$). This $W$ is not semi-invariant for the full torus action on $Y$, but it is semi-invariant for a rank 4 subtorus $T$. Or we can say that $W$ on $\C^6$ is semi-invariant for a rank 5 torus. Choosing an R-charge amounts to choosing six weights $\alpha_s,..., \alpha_z$ such that $W$ has weight $2$. 

\begin{lem}\label{lem.5dflop1} For any decoration $(\alpha, \beta)$ on the graph $\Delta$ we have an equivalence
$$\cC_\Delta \cong D^b_\alpha (Y, W) $$
for some choice of R-charge on $Y$.
\end{lem}
\begin{proof}
 Label $\Delta$ as in Figure \ref{fig.IHmoves}(a). We can cover $Y$ with two toric charts corresponding to $s\neq 0$ or $t \neq  0$, both are isomorphic to $\A^5$. In the first chart the co-ordinates are $y, z$ and
$$x=sp, \quad \hat{y}=sq, \quad u = t/s $$
and the superpotential is $ W = \hat{y} z + uxy$. By Kn\"orrer periodicity (Remark \ref{rem.KP}(4)) the category of matrix factorizations on this chart is equivalent to $D^b_\alpha(\A^3, uxy)$. Similarly matrix factorizations on the other chart are equivalent to $D^b_\alpha(\A^3, vzw)$ where $v=s/t$ and $w=qt$. On the overlap we get the localizations of these categories at $u$ and $v$. It follows that $D^b_\alpha(Y,W)\cong \cC_\Delta$ for some choice of decoration on $\Delta$. 

Moreover it's easy to check that we can pick our R-charge so that it matches any given decoration on $\Delta$; the three equations to solve are
$$\alpha_s+ \alpha_p = \alpha_x, \quad\quad \alpha_q + \alpha_t = \alpha_q, \aand \alpha_t - \alpha_s = \alpha_u$$
and the $\beta$ data is irrelevant. 
\end{proof}

We first introduced the graph $\Delta$ in the context of IH moves, where we replace $\Delta$ with the graph $\Delta'$ from Figure \ref{fig.IHmoves}(b). The algebro-geometric analogue of this is of course a standard flop. The threefold $X$ has another birational model $X'$, it's the other GIT quotient where the unstable locus is $\{p=q=0\}$. We can cross with $\A^2$ and get a family of standard flops changing $Y$ into $Y'=X'\times \A^2$. Looking at Example \ref{eg.3foldflop} or the proof of Lemma \ref{lem.5dflop1} it's evidently natural to identify the category $\cC_{\Delta'}$ with $D^b_\alpha(X', W)$ (if condition \eqref{eq.flopcondition} holds) or $D^b_\alpha(Y', W)$.

\begin{rem}\label{rem.alphatoalpha'} If we choose $\alpha$ data for $\Delta$ then there is a corresponding $\alpha'$ data for $\Delta'$, since both are determined by their values on the external edges.  Obviously condition \eqref{eq.flopcondition} is preserved under this correspondence. Recall that we need this condition to lift $\alpha$ to an R-charge on $\C^4$ and hence on $X$ and $X'$. 
\end{rem}

 Obviously the two graphs $\Delta$ and $\Delta'$ are isomorphic but the point of an IH move is that we have a fixed bijection between $\partial \Delta$ and $\partial \Delta'$ and the isomorphism doesn't respect this. Hence when we embed this move in a larger graph it may change the isomorphism class. Similarly $X$ and $X'$ are isomorphic but there isn't an isomorphism which restricts to the identity away from the flopping curves. But the equivalence
\beq{eq.flopeq}D^b_\alpha(X, W) \isoto D^b_\alpha(X', W) \eeq
provided by Theorem \ref{thm.VGIT} does have this property; see Remark \ref{rem.VGIT}(1). If we consider two larger 3-folds related by a standard flop then they are derived equivalent but they may not be isomorphic. 

This suggests that the categories $\cC_{\Delta}$ and $\cC_{\Delta'}$ should be equivalent relative to their external edges, \emph{i.e.} that we can find an equivalence which intertwines with the restriction functors at each external edge (see Section \ref{sec.categoriesfromgraphs}). This statement will be a key ingredient in Section \ref{sec.invariance}. 

If our $\alpha$ data satisfies \eqref{eq.flopcondition} then we can use the identifications $\cC_\Delta\cong D^b_\alpha(X,W)$ and $\cC_{\Delta'}\cong D^b_\alpha(X', W)$ and the result follows more-or-less immediately from \eqref{eq.flopeq}. Localizing at any external edge deletes both flopping curves and then $\Phi$ becomes the identity.

It remains to prove the statement without condition \eqref{eq.flopcondition} using our 5-fold models $Y$ and $Y'$. These models need some additional uses of Kn\"orrer periodicity and this forces us to confront another issue which we've been ignoring: the $\beta$ data.

 Without a choice of $\beta$ and $\beta'$ data on our graphs the functors at external edges are not quite defined, so the best we could say is that $\cC_\Delta$ and $\cC_{\Delta'}$ are `equivalent relative to their external edges up to shifts'. This is essentially saying that they are equivalent over the ring $R$ from Remark \ref{rem.externalring}. Or we can express it as follows:

\begin{prop}\label{prop.5dflop2} Let $\Delta$ and $\Delta'$ be the two trees from Figure \ref{fig.IHmoves}. Fix a decoration $(\alpha, \beta)$ on $\Gamma$ and let $\alpha'$ be the corresponding data on $\Gamma'$. Then for some choice of $\beta'$ we have an equivalence
$$\Psi: \cC_{\Delta} \isoto \cC_{\Delta'} $$
relative to the external edges. 
\end{prop}
The problem of determining $\beta'$ in terms of $\beta$ is essential for understanding how the category $\cC_\Gamma$ depends on $\beta$ for high genus graphs, but we defer this problem to Section \ref{sec.invariance}. 

\begin{proof}
We can lift both $\alpha$ and $\alpha'$ to the same R-charge $\alpha$ on $\C^6$. Then from Theorem \ref{thm.VGIT} we have an equivalence:
$$\Phi: D^b_\alpha(Y, W) \isoto D^b_{\alpha}(Y', W)$$
By Lemma \ref{lem.5dflop1} the first is equivalent to $\cC_\Delta$ and the second to $\cC_{\Delta'}$. 

Now consider the external edge $x$. Since $\Phi$ is linear over the ring of invariants on $\C^6$ it restricts to give an equivalence between on the open subsets where $x=sp \neq 0$:
$$\Phi: D^b_\alpha(Y_{sp\neq 0}, W) \isoto D^b_{\alpha}(Y'_{sp\neq 0}, W)$$
Our equivalence $D^b(Y, W)\cong \cC_\Delta$ restricts (by construction) to a Kn\"orrer equivalence between $D^b_\alpha(Y_{sp\neq 0}, W)$ and $\cC_x = D^b_\alpha(\C[x^{\pm 1}])$. The same is true on the $Y'$ side. Hence $\Phi$ does intertwine with some equivalence $\cC_x \isoto \cC_x$. Moreover both $\Phi$ and the Kn\"orrer equivalences are linear over $x=sp$, and the only autoequivalences of $\cC_x$ which are linear over $x$ are shifts.

The edge $w$ is identical but the $y$ or $z$ edges need one extra step. Over $\{y\neq 0\}$ we still have a family of standard flops $Y_{y\neq 0} \dashrightarrow Y'_{y\neq 0}$, and $\Phi$ does restrict to give an equivalence here. The functor from $D^b_\alpha(Y, W)$ to $\cC_y$ can be factored as two restriction maps followed by a Kn\"orrer equivalence:
$$D^b_\alpha(Y, W) \to D^b_\alpha(Y_{y\neq 0}, W)\to D^b_\alpha(Y_{ys\neq 0}, W) \isoto \cC_y$$
The second restriction map is an equivalence because there no critical points in the subset $\{y\neq 0, s=0\}$, and it is linear over $y$. We can do a similar factorization on the $Y'$ side and it follows that $\Phi$ intertwines with a $y$-linear equivalence $\cC_y \isoto \cC_y$. 
\end{proof}

\subsection{Orbifold nodes}\label{sec.orbifoldnodes}

Another kind of one-dimensional mirrors are given by nodes with cyclic orbifold structure, \emph{i.e.} substacks of the form
$$ \cZ_{n,k} = \{xy=0\} \subset [\A^2 / \Z_n] $$
where the $\Z_n$ acts diagonally with weights $(1,k)$. These appear in work of Polishchuk and Lekili \cite{LPnodes}. The cases where $k=n-1$ are called `balanced'; they were considered earlier \cite{STZ}. We can glue these orbifold nodes together into chains or rings as in Section \ref{sec.STZ}. 

We don't have much to say about this construction, we're just going to give a brief algebro-geometric argument which connects them to Hori-Vafa mirrors. 
\pgap

By equivariant Kn\"orrer periodicity (Remark \ref{rem.KP}(2)) we have an equivalence:
$$D^b_\alpha(\cZ_{n,k}) \isoto D^b_\alpha\!\big([\A^3 / \Z_n], xyz\big) $$
Here $\Z_n$ must act on $z$ with weight $n-1-k$ so that $W=xyz$ is invariant. We have also chosen an R-charge $\alpha$, the standard choice would be to set $\alpha_x=\alpha_y=0$ and $\alpha_z=2$ and hence recover the ordinary derived category of $\cZ_{n,k}$. 

This 3-dimensional orbifold is `Calabi-Yau' in that $\Z_n$ is acting with trivial determinant. In this sense it is already a Hori-Vafa mirror, if we allow the latter to include toric orbifolds as well as varieties.

 We can also use VGIT to connect it to a variety. The toric polytope for this orbifold is a triangle, typically with lattice points in the interior and on the edges. If we choose a complete triangulation of this lattice polytope then we get the toric data for some Calabi-Yau 3-fold variety $X$. Moreover, $X$ and $[\A^3 / \Z_n]$ are related by VGIT, they are two different (generic) GIT quotients of some vector space by a torus. So by Theorem \ref{thm.VGIT} we have an equivalence:
$$D^b_\alpha\!\big([\A^3 / \Z_n], xyz\big)  \cong D^b_\alpha(X, W) $$
We saw an example of this in Section \ref{sec.3pt}, with $n=3$ and $X=K_{\P^2}$. 

In the balanced case the variety $X$ we end up with is $\A^1_z$ crossed with the small resolution of the $A_{n-1}$ singularity, this is just the McKay correspondence crossed with a line. It's easy to see that the graph $\Gamma$ in this case has genus zero so $\cZ_{n, n-1}$ is mirror to a punctured sphere. In the non-balanced case the dependence of the genus on $n$ and $k$ is more complicated.
\pgap

From our perspective this description of the mirror is rather limited because it only works for a very restricted set of decorations on $\Gamma$, or line fields on the mirror punctured surface. We only have a $\Z^2$ of possible $\alpha$ data, corresponding to the R-charges of $x$ and $y$, and our $\beta$ data must be the canonical one described in Section \ref{sec.HoriVafa}. Perhaps for other choices of $\beta$ we can interpret $\cC_\Gamma$ as a twisted version of $D^b_\alpha(\cZ_{n,k})$ but we do not know how to make sense of this.

\subsection{Nodal elliptic curves} \label{sec.nodalelliptic}

Let $\Sigma$ be a torus with a single puncture. A pair-of-pants decomposition of $\Sigma$ produces a trivalent graph $\Gamma$ of genus one with a single external edge, so $\Gamma$ must have one vertex and a loop. 

Let us label the half-edges of $\Gamma$ by $x,y,z$, where $z$ is the external edge and $x\cup y$ the internal edge. Now add a decoration. We are forced to set $\alpha_z=2$ and 
$$\alpha_x = -\alpha_y = a $$
for some $a\in \Z$. Then up to trivial modifications the $\beta$ data is determined by the invariant:
$$ \beta_{yx} - \beta_{xy} =  \gamma_{yx} =   \delta_{xy} \in \Z/a$$
This is a special case of the situation considered in Section \ref{sec.STZ} and as we did there we can construct a one-dimensional mirror curve $C$. In this case $C$ is formed taking a single node and then gluing it to itself. The result is a nodal elliptic curve, \emph{i.e.} the unique plane curve $C$ of arithmetic genus 1 with a single node. $C$ can also be described as a rational curve with two points glued together.\footnote{If you 
prefer to call $C$ a rational curve instead of an elliptic curve then feel free to do so.}

There is a rank 1 torus acting on $C$ and the $\alpha$ data specifies an R-charge $a: \C^* \to \C^*$.  If we choose an integer $b$ representing $\delta_{xy}$ then we have a class in  $H^1_{\acute{e}t}(C, \Z)= \Z$. Then as in Section \ref{sec.STZ} we have an equivalence
\beq{eq.predictedFM}\cC_{\Gamma} \cong D^b_a(C, b) \eeq
to the twisted derived category of $C$. The choices $a=b=0$ give the usual derived category $D^b(C)$.  In this sense $C$ is the mirror to $\Sigma$. 
\pgap

The category $D^b_a(C, b)$ obviously only depends on the residue $[b]\in \Z/a$, see Remark \ref{rem.bmodalpha}. But there another more subtle equivalence here which only becomes obvious on the mirror side. 

line fields on the surface $\Sigma$ are specified by their winding numbers along two loops generating $H_1(\Sigma)$, these are of course the numbers $a$ and $b$. The Fukaya category is obviously preserved under the action of the mapping class group $SL_2(\Z)$ so it will not change if we replace $(a,b)$ with $(b, -a)$. Hence if we had actually proved mirror symmetry we could immediately deduce that we have an equivalence
$$D^b_a(C, b) \isoto D^b_b(C, -a)$$
exchanging the roles of the R-charge and the $H^1$ twist. Instead we will give a proof that this equivalence does indeed hold, but we'll do it using only B-side arguments.

In the context of smooth elliptic curves this equivalence is very familiar, it's the original Fourier-Mukai transform between an elliptic curve and its Jacobian. Here it is an autoequivalence, it exchanges rank and degree, and is well-known to be mirror to this element of the mapping class group.  Burban and Kreu{\ss}ler have shown that this extends to the nodal elliptic curve \cite{BK}; the Jacobian of $C$ is isomorphic to $C$ and the universal sheaf gives a derived equivalence. So our task is simply to incorporate the R-charge and twisting.

\pgap

The curve $C$ has the following \'etale covering: 
\beq{etalecover}\begin{tikzcd}[column sep = huge]
\A^1_t\setminus 0  \arrow[r,  shift left=1ex, "i:\, t\mapsto {(t,0)}"] \arrow[r, shift right=1ex, "j:\, t\mapsto {(0,t^{-1})}"'] & Z = \operatorname{Spec} \C[x,y]/xy
\end{tikzcd}\eeq
So a sheaf on $C$ can be described as a sheaf $\cE$ on the node $Z$ together with an isomorphism between $i^*\cE$ and $j^*\cE$. An object of the twisted category $D^b_a(C, b)$ can be described as an object in $D^b_a(Z)$,  \emph{i.e.}~a graded module $\cE$ with a differential,  together with a quasi-isomorphism:

$$\psi: i^*\cE \isoto j^*\cE[b]$$

\begin{eg}\label{eg.objectsonC}\hspace{1pt}

\begin{enumerate}
\item Suppose $b = ka$ for some $k$. Then $D^b_a(C, b)$ contains a family of `line-bundles' which are locally the structure sheaf, but have transition functions
$$ \lambda t^k : \cO \isoto \cO[b] $$
for some $\lambda\in \C^*$. 

\item  A torsion sheaf defines an object of $D^b_a(C, b)$ for any $b$, but if $a\neq 0$ then all torsion sheaves must be supported at the node. 

\item For any $a, b$ we have an object given by $\cE = \cO_{x=0}\oplus \cO_{y=0}[-b]$, with $\psi$ the identity map. This is a twisted version of the push-down of the structure sheaf of the normalization.
\end{enumerate}
The objects (1) are mirror to loops in $\Sigma$ along which the line field has zero winding number. The same is true for torsion sheaves at smooth points when $a=0$. Torsion sheaves at the node are mirror to arcs, as is the object (3). 
\end{eg}

We will produce the equivalence \eqref{eq.predictedFM} by writing down a Fourier-Mukai kernel for it. On $C\times C$ the possible R-charges are specified by two integers $(a_1, a_2)$, and the twisting class by another two integers $(b_1, b_2)$. An object in  $D^b_{a_1, a_2}(C\times C, b_1, b_2)$ produces a functor:
$$D^b_{a_1}(C,  -b_1) \To D^b_{a_2}(C,  b_2)$$

We can explicitly describe such an object by using the cover \eqref{etalecover} for both factors of $C$. So it's an object $\cE$ on the affine variety $Z\times Z$, together with quasi-isomorphisms
$$\psi_1: (i\times \id)^*\cE \isoto (j\times \id)^*\cE [b_1] \ \aand \psi_2: (\id\times i)^*\cE \isoto (\id\times j)^*\cE[b_2] $$
satifying the evident compatibility relation. In fact our kernel will just be a sheaf (with no differential) and $\psi_1, \psi_2$ will be isomorphisms.

Now let $(a_1, a_2)$ be any two integers and let $\cE$ be the sheaf on $Z\times Z$ which is the cokernel of the following map:
$$\begin{tikzcd}[column sep = huge, ampersand replacement=\&]
\cO\oplus \cO[a_2-a_1]  \arrow[r,  "{\begin{pmatrix} y_1 & x_2 \\ y_2 & x_1 \end{pmatrix}}"]  \& \cO[-a_1] \oplus \cO[a_2 ]
\end{tikzcd}$$
If $a_1=a_2$ then $\cE$ is the ideal sheaf of the diagonal, but if they are not equal then the diagonal is not invariant and this is meaningless.

Away from the origin $\cE$ is a line-bundle so locally it's a shift of $\cO$. But the complement of the origin is not connected, and the shift is different on the different components. For example if we pull-back along $i\times\id $ then $\cE$ becomes isomorphic to  $\cO$, using the exact sequence:
$$\begin{tikzcd}[column sep = large, ampersand replacement=\&]
\cO\oplus \cO[a_2-a_1]  \arrow[r,  "{\begin{pmatrix} 0 & x_2 \\ y_2 & t_1 \end{pmatrix}}"]  \& \cO[-a_1] \oplus \cO[a_2 ]  \arrow[r,  "{\begin{pmatrix} t_1 \\ -y_2  \end{pmatrix}}"] \& \cO \arrow[r] \& 0
\end{tikzcd}$$

But if we pull-back along $j\times \id$ then we can use a similar sequence to identify $\cE$ with $\cO[a_2]$.  For our first transition function we take the identity map on $\cO$ and view it as an isomorphism:
$$\psi_1:  (i\times \id)^*\cE \stackrel{1}{\longrightarrow} (j\times \id)^* \cE [ -a_2] $$

Our second transition function is defined similarly, we identify $(\id\times i)^*\cE$ with $\cO[-a_1]$ and $(\id\times j)^*\cE$ with $\cO$, and set:
$$\psi_2 : \; (\id\times i)^*\cE  \;\stackrel{1}{\longrightarrow} \; (\id\times j)^*\cE [-a_1] $$

It's obvious that $\psi_1$ and $\psi_2$ commute with each other (after further restriction) so we have defined an object
$$E \;\in D^b_{a_1, a_2}(C\times C, -a_2, -a_1) $$
and hence a functor from $D^b_a(C, b)$ to $D^b_b(C, -a)$ for any $a, b$. 

\begin{prop}\label{prop.FM} The Fourier-Mukai kernel $E$ defined above gives an equivalence:
$$\mathfrak{E}: D^b_a(C, b) \isoto D^b_b(C, -a) $$
\end{prop}
\begin{proof}
We have our kernel $E\in D^b_{a,b}(C\times C, -b, -a)$. Dualising and swapping the two factors gives an object $E^\vee \in D^b_{b,a}(C\times C, a, b)$. We claim that the inverse to $E$ is $E^\vee[1]$. 

The kernel for the identity functor on $D^b_a(C, b)$ is the object in $D^b_{a, a}(C\times C, -b, b)$ given by 
$$ \mathcal{I}_b =  \frac{\cO}{(x_1-x_2, y_1-y_2)} \;\oplus \; \frac{\cO}{x_1y_2-1}[b]  \;\oplus \; \frac{\cO}{y_1x_2-1}[-b] $$
with the obvious transition functions. It's easy to verify that:\footnote{This is essentially verifying that that Serre duality holds in this twisted derived category.}
$$ (\cI_b)^\vee = \cI_{-b}[-1] \quad \in D^b_{a,a}(C\times C, b, -b) $$
In the category $D^b_{a,b,a}(C\times C\times C, -b, 0, b)$ there is an evaluation map
$$\pi^*_{12}E\otimes\pi_{23} E^\vee \To \pi^*_{13}(\cI_b) $$
where $\pi_{ij}$ the projection to the $i$th and $j$th factors. Taking duals and applying adjunction yields a map
$$ \cI_{-b} \To E^\vee\star E[1]$$
 where $\star$ denotes convolution. In the case $a=b=0$ the kernel $E^\vee[1]$ is the inverse to $E$ and this map is an isomorphism \cite{BK}, hence it is an isomorphism for any $a,b$. An identical argument produces an isomorphism $\cI_b\isoto E\star E^\vee[1]$.
\end{proof}

\begin{cor} \label{cor.SL2} The category $D^b_a(C, b)$ only depends on $\gcd(a,b)$. 
\end{cor}
\begin{proof} The transformations $(a,b) \mapsto (a, b+a)$ and $(a,b)\mapsto (b, -a)$ generate $SL_2(\Z)$. 
\end{proof}

This is our algebro-geometric version of the result \cite{LPinvs} that a line field on a punctured torus is determined, up to the action of the mapping class group, by the gcd of the winding numbers. 

In Section \ref{sec.genus0or1} we're going to extend this result to general graphs of genus one. To do this we're going to need a technical lemma about the equivalence $\mathfrak{E}$.
\pgap

There's a restriction functor from $D^b_a(C, b)$ to $D^b_a(Z)$, and hence from there to $D_{sg, a}(Z)$. This latter category is equivalent to $D^b(\C[z^{\pm 1}]) $ where $z$ has degree 2 which is the 2-periodic category of vector spaces (Remark \ref{rem.localizationfornodes}). Clearly it would be the same thing as $D_{sg, a}(C, b)$ if we knew the definition of a twisted derived category of singularities.

We can ask if the equivalence $\mathfrak{E}$ intertwines with a functor $D^b_a(Z)\to D^b_b(Z)$. This is obviously not true, but it becomes true if we pass to $D_{sg}$. 

\begin{lem}\label{lem.FMcommuteswithDsg} The equivalence $\mathfrak{E}$ from Proposition \ref{prop.FM} intertwines with an equivalence $D_{sg, a}(Z) \isoto D_{sg, b}(Z)$.
\end{lem}
Note that since both of these categories are equivalent to $D^b_{\Z_2}(Vect)$ there are only two possible equivalences between them, differing by a shift. It's not significant which of these equivalences occurs,  we could always compose $\mathfrak{E}$ with a shift. 

The category $D^b_a(C, b)$ carries an action of $\C[z]$ and passing to $D_{sg, a}(Z)$ is localizing the category away from $z=0$. So really the statement we should prove is that $\mathfrak{E}$ is linear over $z$, indeed there is a shadow of this statement in the proof below. Alternatively if one could prove that $D_{sg, a}(Z)$ was the quotient of $D^b_a(C,b)$ by the subcategory of compact objects then the lemma would follow simply because $\mathfrak{E}$ is an equivalence.

\begin{proof} Consider the object $\cO_0\in D^b_a(C,b)$ (Example \ref{eg.objectsonC}(2)). Mapping this to $D_{sg, a}(Z)$ produces the $\Z_2$-graded vector space $V=\C\oplus \C[1]$ which is a split-generator.  The endomorphism algebra of $\cO_0$ is
$$\End(\cO_0) = A =\C\langle\theta_1, \theta_2\rangle/(\theta_1^2, \theta_2^2)$$
where $|\theta_1|=1-a$ and $|\theta_2|=1+a$. We know that $A$ must be a $\C[z]$-algebra, and since there's a unique central element of degree 2 we must have $z=\lambda(\theta_1\theta_2+\theta_2\theta_1)$ for some $\lambda\in \C$. Moreover $\lambda\neq 0$ since $\cO_0$ doesn't map to zero in $D_{sg, a}(Z)$.    If we invert $z$ then $A$ becomes $\End(V)$ which up to grading is a rank two matrix algebra. The functor $D^b_a(C, b)\to D_{sg, a}(Z)$ can also be described as
$$\C[z^{\pm 1}]\otimes_{\C[z]}\Hom(\cO_0, -): \; D^b_a(C, b)\To D^b(A_{z^{-1}})$$
followed by the Morita equivalence between $A_{z^{-1}}$ and $\C[z^{\pm 1}]$. 

Now consider the object $\mathfrak{E}(\cO_0)\in D^b_b(C, -a)$. It is easy to compute that this is the sheaf $\cO/y[-a]\oplus \cO/x[b]$ with the transition function $t$ (this is a combination of Examples \ref{eg.objectsonC}(1) and (3)). 

Let $A' = \End(\mathfrak{E}(\cO_0))$. We know that $A'$ is isomorphic to $A$ as an algebra  because $\mathfrak{E}$ is an equivalence. We do not know \emph{a priori} that $A'\cong A$ as $\C[z]$ algebras, but from the structure of $A$ the $z$ actions must agree up to a scalar, and since $\mathfrak{E}(\cO_0)$ evidently does not map to zero in $D_{sg, b}(Z)$ the scalar is non-zero.  Hence $\mathfrak{E}$ induces a well-defined isomorphism 
$$\mathfrak{E}: A_{z^{-1}} \isoto A'_{z^{-1}} $$
on the localized algebras and the lemma follows.
\end{proof}

\begin{rem}\label{rem.compacttorus} It would be nice to extend our constructions to a compact torus but we're not sure how to mirror the line fields. The mirror of a compact torus is a smooth elliptic curve, or perhaps more precisely the formal family of elliptic curves degenerating to the cuspidal one. It seems plausible that R-charge could be replaced somehow with the action of an elliptic curve on itself. We do not know what the analogue of the class $b$ is. 
\end{rem}

\section{$\cC_\Gamma$ as an invariant of decorated graphs}\label{sec.invariance}

In Section \ref{sec.categoriesfromgraphs} we explained how to construct a dg-category $\cC_\Gamma$ from a decorated trivalent graph $\Gamma$. We also explained that without a decoration there is still a canonical $\Z_2$-graded dg-category $\cD_\Gamma$. This latter category is equivalent to the Fukaya category of the surface $\Sigma$ corresponding to $\Gamma$, so it follows that the category only depends on the genus and number of external edges of $\Gamma$. In this section we'll give an alternative proof of this fact using only the algebraic geometry of the previous sections, plus some elementary combinatorics.

 In the $\Z$-graded case the Fukaya category also depends on the line field, or more precisely on the orbit of the line field under the mapping class group of $\Sigma$. In \cite{LPinvs} the authors give a complete description of the numerical invariants of line fields. We will fshow how to extract these same invariants combinatorially from a decorated graph, then we prove that the category $\cC_\Gamma$ only depends on the topology of $\Gamma$ plus these invariants. 
\pgap

Before we start we recall the following elementary and well-known lemma, and also its proof since we will need a detail from it.\footnote{We thank Dan Petersen for telling us this argument.}

\begin{lem}\label{lem.IHmoves} Let $\Gamma$ and $\Gamma'$ be two connected trivalent graphs of the same genus, and suppose we have a chosen bijection between $\partial\Gamma$ and $\partial \Gamma'$. Then there is a sequence of IH moves that transforms $\Gamma$ into $\Gamma'$. 
\end{lem}

\begin{proof} Firstly suppose both graphs are trees and pick an ordering $1,..., n$ of the external edges. Now consider a third `straight tree'  $\Gamma''$ (as in Figure \ref{fig.straighttree}) which has this same set of external edges, such that (a) external edges 1 and 2 meet at a vertex, as do external edges $n$ and $n-1$, and (b) for $2\leq j \leq n-2$ the external edges $j$ and $j+1$ are connected by a single  internal edge. We claim that $\Gamma$ can be transformed into $\Gamma''$ by a sequence of IH moves, hence $\Gamma$ can also be transformed into $\Gamma'$. 

Turning $\Gamma$ into $\Gamma''$ is a simple recursive procedure. Consider the path of internal edges joining external edges 1 and 2. If this path is non-zero  then pick one of the edges in it and perform an IH move to make the path shorter. Repeat until external edges 1 and 2 meet at a vertex, and denote the internal edge leaving that vertex by $2'$. Use the same procedure to make $2'$ meet the external edge 3 at a vertex. Repeat.

Finally suppose $\Gamma$ and $\Gamma'$ have genus bigger than zero. On each graph we can cut internal edges until we get trees, this creates more external edges. Extend the given bijection to the additional external edges arbitrarily then apply the result for trees.
\end{proof}

\subsection{Invariance in the $\Z_2$-graded case}\label{sec.Z2invariance}

As we saw in Section \ref{sec.Z2gradedcase} on any trivalent graph there is a canonical choice of weak decoration so there is an associated $\Z_2$-graded category. We will now give a purely algebro-geometric proof that this category only depends on the genus $g$ and number $n$ of external edges of the graph. 

In outline the proof is very simple: any two graphs with the same $g$ and $n$ can be connected to each other by a sequence of IH moves, and the derived equivalence for the standard flop (Section \ref{sec.flops}) shows that an IH move doesn't change the category.  However it is slightly more subtle than it appears. There are other choices of weak decoration and they give (presumably) different categories, so we must take care when performing an IH move or flop to keep track of the weak decoration. Doing this in detail also gives us a warm-up for the $\Z$-graded case which we consider in the following sections.
\pgap

Let $\Delta$ be a genus zero trivalent graph with four external edges, the relevant graph for an IH move. We're interested in weak decorations on $\Delta$ modulo trivial modifications. If we allow all trivial modifications then by Lemma \ref{lem.H1mod2} there is only one equivalence class. Instead let's only allow relative trivial modifications, \emph{i.e.}~modifications of types (V) and (I), but not (E). Let's write $\cG_\Delta$ for the resulting set of equivalence classes. If we choose an element $b\in \cG_\Delta$ then we can construct a category $\cD_\Delta$ together with functors 
$$\cD_\Delta \To \cD_x $$
at each of the four external edges. If we change $b$ then we get an equivalent category but the functors at the external edges change by shifts. 

The set $\cG_\Delta$ is obviously a torsor for the relative cohomology group:
$$H^1(\Delta, \partial\Delta; \Z_2)\; \cong (\Z_2)^{\oplus 3}$$
Unlike the absolute case this torsor is not canonically trivial; the proof of Lemma \ref{lem.H1mod2} fails because the invariants $b_c$ are not additive on relative cycles. But we can still use these invariants to record the class of our weak decoration in the following way. 

Let us label the half-edges of $\Delta$ as in the left-hand-side of Figure \ref{fig.IHmoves}. Between any pair of external edges there is a path, involving either 2 or 4 half-edges, and we can add the $\beta$'s along each path to get six elements of $\Z_2$:
$$\epsilon_{xy} = \beta_{xy} + \beta_{yx}, \quad   \epsilon_{xz} = \beta_{xu} + \beta_{ux} + \beta_{vz} + \beta_{zv}, \quad ... $$
These $\epsilon$'s are obviously invariant under relative trivial modifications and  obey relations of the form:
$$\epsilon_{xy} + \epsilon_{xz} + \epsilon_{yz} = 1$$
They identify $\cG_\Delta$ with a 3-dimensional affine subspace of $(\Z_2)^6$. 

Now perform an IH move, replacing $\Delta$ with a different graph $\Delta'$ having the same set of external edges. The $\epsilon$ invariants provide a canonical bijection $\cG_{\Delta}\isoto \cG_{\Delta'}$, so we can declare that a weak decorations on $\Delta$ and $\Delta'$ are \emph{equivalent} if their $\epsilon$ invariants agree. 

\begin{prop}\label{prop.mod2flop} Choose a weak decoration $\beta$ on $\Delta$ and an equivalent weak decoration $\beta'$ on $\Delta'$. Then the categories $\cD_\Delta$ and $\cD_{\Delta'}$ are equivalent relative to their external edges.
\end{prop}
\begin{proof}
This equivalence is the $\Z_2$-graded version of the 3-fold flop, see Section \ref{sec.flops}. The delicate part is verifying that the flop equivalence intertwines the functors at external edges.

Using trivial modifications of type (V) we can arrange that 
$$\delta_{xz}= \beta_{ux}+\beta_{vz}=0 \aand \delta'_{xy} = \beta'_{u'x}+\beta'_{v'y}=0$$
and $\beta'_{xz}=\beta_{xu}$. Then since $\beta$ and $\beta'$ are equivalent it follows that:
$$\beta'_{yv'} = \beta_{yx}, \quad\quad \beta'_{zx} = \beta_{zv}, \quad\quad \beta'_{wy} =\beta_{wv} $$
It's clear that we have an equivalence $K: D^b(X, W)\isoto \cD_\Gamma$. Moreover since $\delta_{xz}=0$ we can choose $K$ such that on both toric charts $K$ is the identity functor, rather than a shift (\emph{c.f.}~the proof of Lemma \ref{lem.HVmirror}). Similarly we have an equivalence $K': D^b(X', W) \isoto \cD_{\Gamma'}$ which is the identity on both charts. We have an equivalence $\Phi: D^b(X, W)\isoto D^b(X', W)$, so we just need to check that $\Psi = K'\circ \Phi\circ K^{-1}$ commutes with the functors at each external edge. 

Consider the edge $x$. We have canonical isomorphisms
$$ \C^2_{y, u}\times \C^*_x \; \cong \; X_{x\neq 0} \; \cong \; X'_{x\neq 0} \; \cong \; \C^2_{z, u'} \times \C^*_x $$
under which $y\mapsto u'x$ and $u\mapsto zx^{-1}$. Since $\Phi$ is the identity away from the exceptional curve, and $K$ and $K'$ are the identity functor on each chart, we have (after inverting $x$) that:
$$K'\circ \Phi\circ K^{-1}: \cO/y \mapsto \cO/u' $$
Applying the edge functor $\cC_\Gamma\to D^b([x^{\pm}])$ maps this object to $\cO[\beta_{xu}]$, and doing the same on the $\Gamma'$ side produces $\cO[\beta'_{xz}]$. So the functors commute. The argument at the remaining external edges is identical.
\end{proof}

Now take any connected trivalent graph $\Gamma$ and perform an IH move to get a new graph $\Gamma'$. There's a canonical isomorphism $H^1(\Gamma, \Z_2) \cong H^1(\Gamma', \Z_2)$ so by Lemma \ref{lem.H1mod2} for any class of weak decorations on $\Gamma$ there's a corresponding class of weak decorations on $\Gamma'$.

Suppose we choose an explicit weak decoration $\beta$ on $\Gamma$, and let $\Delta\subset \Gamma$ and $\Delta'\subset \Gamma'$ be the half-edges involved in the IH move. We can get a weak decoration $\beta'$ on $\Gamma'$ by setting $\beta'=\beta$ where the two graphs are identical, and on $\Delta'$ choosing some $\beta'$ which is equivalent to $\beta$ in the sense above. It's obvious that this prescription is well-defined with respect to trivial modifications and the classes of $\beta$ and $\beta'$ in $H^1$ will agree. In particular if $\beta$ represents the canonical class of weak decorations on $\Gamma$ then the same is true for $\beta'$ on $\Gamma'$. 

\begin{thm}The category  $\cD_\Gamma$ only depends on the genus and number of external edges of $\Gamma$.
\end{thm}
\begin{proof} 
Let $\Gamma'$ be a second trivalent graph with the same genus and number of external edges and fix a bijection between $\partial\Gamma$ and $\partial\Gamma'$. By Lemma \ref{lem.IHmoves} there's a sequence of IH moves connecting the two graphs. Choose a weak decoration on $\Gamma$ representing the zero class in $H^1(\Gamma, \Z_2)$, and for each IH move modify the weak decoration as described above; the end result is a weak decoration on $\Gamma'$ representing the zero class. 

Using the weak decoration we can construct the category $\cD_{\Gamma}$ as a homotopy limit over a diagram, and after each IH move we get a new category. But when we perform an IH move only part of the diagram changes and by Proposition \ref{prop.mod2flop} the limit over that part of the diagram does not change; hence the overall limit doesn't change. 
\end{proof}

\begin{rem}One could also choose a non-zero class $b\in H^1(\Gamma, \Z_2)$ and consider the corresponding $\Z_2$-graded category $\cD_{\Gamma,b}$. By the proof above this category is invariant under IH moves. Moreover, it's not hard to show that for any $b'\neq 0 \in H^1(\Gamma', \Z_2)$ there's a sequence of IH moves turning $b'$ into $b$, which means that $\cD_{\Gamma, b}$ also only depends on the topology of $\Gamma$, and not the specific class $b$ (as long as $b\neq 0$). 

For a proper discussion here we should consider all $\Z_2$ decorations, which means allowing the $\alpha$ data to be non-zero as well. We won't bother to do this since the $\Z$-graded case is more interesting. 
\end{rem}

\subsection{Decorations and IH moves}

We now want to upgrade the result of the previous section to the $\Z$-graded case. The first step is to understand how decorations should change under IH moves.

Once again let $\Delta$ and $\Delta'$ be the two graphs from Figure \ref{fig.IHmoves}. If we have a decoration $(\alpha, \beta)$ on $\Delta$ then we need to find an `equivalent' decoration $(\alpha', \beta')$ on $\Delta'$.  The change in $\alpha$ is easy and was already discussed in Remark \ref{rem.alphatoalpha'}: both $\alpha$ and $\alpha'$ are determined by their values on the external edges so we just need to insist that they agree on $x,y,z,w$. Since $\alpha'$ is determined by $\alpha$ in such a simple way we will drop the notation $\alpha'$ and write $\alpha$ for the data on both graphs. 

To understand the change in $\beta$ we first we need to understand the orbits of $\beta$ under relative trivial modifications. By definition the set of decorations on $\Delta$ modulo relative trivial modifications form a torsor for the group:

$$G_\Delta^\alpha=  \frac{ \Z_{\alpha_x}\oplus \Z_{\alpha_y}\oplus\Z_{\alpha_z}\oplus\Z_{\alpha_w}\oplus (\Z_{\alpha_u})^{\oplus 2}}{\smat{1 & 1 & 0 & 0 & 1 & 0 \\ 0 & 0 &1 & 1 & 0 & 1 \\ 0 & 0 & 0 & 0 & 1 &-1 }} \;\; \cong \; \;\frac{ \Z_{\alpha_x}\oplus \Z_{\alpha_y}\oplus\Z_{\alpha_z}\oplus\Z_{\alpha_w}}{\smat{1 & 1 & 1 & 1 \\ 0 & 0 &\alpha_u & \alpha_u }} $$

We can choose a trivialization of this torsor as follows. Recall that $\delta_{xw} = \beta_{ux} -\beta_{vw}\in \Z_{\alpha_u}$. From $\beta$ we extract:
\beq{eq.B}B(\beta) = \big(\beta_{xu}, \,\beta_{yx},\, \beta_{zw} + \delta_{xw}, \,\beta_{wv} + \delta_{xw} \big) \; \in G_\Delta^\alpha\eeq
This is a well-defined element of $G_\alpha$, and determines $\beta$ uniquely up to relative trivial modifications. Of course there are many other choices we could have made in defining $B$; this one is convenient for the geometry of the flop.

Now we need a similar bijection for $\beta'$. The orbits of $\beta'$ under relative trivial modifications form a torsor for the group
$$G_{\Delta'}^\alpha=  \frac{ \Z_{\alpha_x}\oplus \Z_{\alpha_y}\oplus\Z_{\alpha_z}\oplus\Z_{\alpha_w}}{\smat{1 & 1 & 1 & 1 \\ 0 & \alpha_{u'} & 0 & \alpha_{u'} }} $$
and we define:
\beq{eq.B'}B'(\beta')  = \big(\beta'_{xu'}, \,\beta'_{yv'}+ \delta'_{xw},\, \beta'_{zu'}, \,\beta'_{wv} + \delta'_{xw} \big) \; \in G_{\Delta'}^\alpha\eeq

\begin{defn}\label{defn.equivalence}\hspace{1cm}

\begin{enumerate}\item[(a)] We declare that the decorations $(\alpha, \beta)$ and $(\alpha, \beta')$ are \emph{equivalent} if $B(\beta)$ and $B'(\beta')$ are equal in the common quotient:
$$G^\alpha =   \frac{ \Z_{\alpha_x}\oplus \Z_{\alpha_y}\oplus\Z_{\alpha_z}\oplus\Z_{\alpha_w}}{\smat{1 & 1 & 1 & 1 \\ 0 & 0 &\alpha_u & \alpha_u\\ 0 & \alpha_{u'} & 0 & \alpha_{u'} }}$$

\item[(b)] Suppose that $\Delta$ is part of a larger graph $\Gamma$ and that $\Gamma'$ is the graph obtained from $\Gamma$ by performing the IH move that turns $\Delta$ into $\Delta'$. We declare that decorations $(\alpha, \beta)$ on $\Gamma$ and $(\alpha, \beta')$ on $\Gamma'$ are \emph{equivalent} if they are equivalent on $\Delta$ and $\Delta'$ and identical elsewhere.
\end{enumerate}
\end{defn}

\begin{rem}\label{rem.coarserER}
Notice that part (b) here is well-defined with respect to trivial modifications: if $\beta$ and $\beta'$ are equivalent and we perform a trivial modification to $\beta$ then there is a corresponding trivial modification for $\beta'$ which restores the equivalence. So it is really an equivalence between orbits of $\beta$ and $\beta'$ under trivial modifications. However it is \emph{not} a bijection: we can have $\beta_1$ and $\beta_2$ on $\Gamma$ which are both equivalent to the same $\beta'$ on $\Gamma'$, without being related to each other by trivial modifications. 

This generates an equivalence relation on decorations which is strictly coarser than the one defined by trivial modifications alone. This is related to spherical twists on the category $\cC_\Gamma$, see Remark \ref{rem.spherical}.
\end{rem}

\begin{rem}\label{rem.epsilon} It would be nice if we could formulate this definition without choosing a trivialization of the torsors on each side, but we haven't found a way to do this. The obvious thing to try is to define (as we did in the $\Z_2$-graded case) an $\epsilon$ invariant for each ordered pair of external edges of $\Delta$, by taking the alternating sum of the $\beta$'s between them. For example we have 
\al{ \epsilon_{wx} := \beta_{wv} - \beta_{vw} + \beta_{ux} - \beta_{xu} &= \gamma_{wv}+\gamma_{ux} \\
&= \beta_{wv} + \delta_{xw} -\beta_{xu} \quad \in \Z/(\alpha_x, \alpha_w, \alpha_u) }
or:
$$\epsilon_{yx}:= \gamma_{yx} = \beta_{yx} - \beta_{xy}\quad \in \Z/(\alpha_x,\alpha_y)$$
These $\epsilon$ residues are antisymmetric and obviously invariant under relative trivial modifications. The six of them together give an element of the group:
$$H_\Delta^\alpha = \frac{\Z}{(\alpha_x, \alpha_y)}\oplus\frac{\Z}{(\alpha_x, \alpha_z, \alpha_u)}\oplus \frac{\Z}{(\alpha_x, \alpha_w, \alpha_u)} \oplus \frac{\Z}{(\alpha_y, \alpha_z, \alpha_u)}\oplus \frac{\Z}{(\alpha_y, \alpha_w, \alpha_u)}\oplus\frac{\Z}{(\alpha_z, \alpha_w)} $$
If we do this on $\Delta'$ we get an element of a different group $H_{\Delta'}^\alpha$, so the $\epsilon$ invariants can only be compared in the common quotient
$$H^\alpha =  \frac{\Z}{(\alpha_x, \alpha_y, \alpha_{u'})}\oplus\frac{\Z}{(\alpha_x, \alpha_z, \alpha_u)}\oplus\; ...\; \oplus\frac{\Z}{(\alpha_z, \alpha_w, \alpha_{u'})} $$
(note that $(\alpha_x, \alpha_w, \alpha_u)=(\alpha_x, \alpha_w, \alpha_{u'})$ and $(\alpha_y, \alpha_z, \alpha_u)= (\alpha_y, \alpha_z, \alpha_{u'})$).
Unfortunately this group is not quite good enough to detect equivalence between $\beta$ and $\beta'$, see Remark \ref{rem.refinedepsilon} below. 
\end{rem}

\begin{rem}\label{rem.bcunderIH} 

Recall that for any cycle $c$ in the larger graph $\Gamma$ we have a residue $b_c \in \Z/I_c$ given by the alternating sum of the $\beta$'s along $c$ \eqref{eq.bc}, and that $b_c$ is invariant under trivial modifications of $\beta$. 

Now perform an IH move, turning $\Gamma$ into $\Gamma'$, and replace $(\alpha, \beta)$ with an equivalent $(\alpha, \beta')$. There are three possibilities for how the cycle $c$ will change.
\begin{itemize}
\item If $c$ avoids the vertices involved in the IH move then $c$ isn't affected and obviously the value of $b_c$ doesn't change.

\item If $c$ includes the edge involved in the IH move then $c$ is transformed into a cycle $c'$ on the new graph $\overline{\Gamma}'$, which is either the same length as $c$ or is one edge shorter. If $c$ and $c'$ are the same length then the ideals $I_c$ and $I_{c'}$ will be equal and it's clear that $b_c=b_{c'}$ since these residues only depend on the $\epsilon$ invariants from the preceeding remark. 

 If $c'$ is shorter than $c$ then $I_c = I_{c'} + (\alpha_{u})$ and:
 $$ b_c = [b_{c'}] \; \in \Z/ I_{c} $$

\item  If $c$ doesn't include the affected edge but does meet one of its vertices then the question is more complicated. We'll return to this in Section \ref{sec.genus2invariance} with some additional hypotheses in place.  
\end{itemize}
\end{rem}

\begin{rem}\label{rem.associatedweak1} Again suppose we peform an IH move turning $\Gamma$ into $\Gamma'$ and we have equivalent decorations  $(\alpha, \beta)$ and $(\alpha, \beta')$ on the two graphs. Assume that all the $\alpha$ values are even; if this holds on $\Gamma$ then it also holds on $\Gamma'$. Then by reducing all the $\beta$ and $\beta'$ values mod 2 we get associated weak decorations on each graph. Since the invariants $\epsilon$ and $\epsilon'$ (from Remark \ref{rem.epsilon}) agree it's clear that the two weak decorations will be equivalent in the sense of Section \ref{sec.Z2invariance}, and in particular they give the same class in $H^1(\Gamma, \Z_2)\equiv H^1(\Gamma', \Z_2)$.
\end{rem}

\begin{rem}\label{rem.refinedepsilon} Suppose we have equivalent decorations $(\alpha, \beta)$ and $(\alpha, \beta')$ on $\Delta$ and $\Delta'$.  Suppose that
$$\alpha_x \equiv \alpha_y \equiv \alpha_z \equiv \alpha_w \equiv 0 \mod 4 $$
so $\alpha_u\equiv \alpha_{u'} \equiv 2 \mod 4$. We've declared that two decorations $\beta$ and $\beta'$ are equivalent if $B(\beta)$ and $B'(\beta')$ are the same in the group $G^\alpha$, which after reducing everything mod 4 becomes:
$$ \frac{ \Z_4^{\oplus 4}}{\smat{1 & 1 & 1 & 1 \\ 0 & 0 &2& 2\\ 0 & 2 & 0 & 2 }} \; \cong \;\Z_4\oplus \Z_2^{\oplus 2} $$
In Remark \ref{rem.epsilon} we considered invariants $\epsilon$ and $\epsilon'$ but these can only be compared in the group $H^\alpha$, which is $\Z_2^{\oplus 6}$ after reducing mod 4. It follows that in this situation the invariants $\epsilon$ and $\epsilon'$ are not strong enough to detect equivalence of decorations.\footnote{In fact this is the only situation where this issure arises; if $\alpha$ is not zero mod 4 at any of the external edges then $H^\alpha$ does detect equivalence. This claim is elementary to prove but we won't need it.} In our formulation this subtlety is  responsible for the appearence of the Arf invariant in genus $\geq 2$.

There is a way to refine the $\epsilon$ invariants here which will be useful in following sections. First, express $\epsilon$ in terms of $B(\beta)$. If we write $(B_x, B_y, B_z, B_w)$ for the components of $B(\beta)$ \eqref{eq.B} then:
\beq{eq.epsilon}\epsilon = \big(B_y - B_x + 1, \; B_z - B_x, \; B_w - B_x, \; B_z-B_y, \; B_w- B_y, \; B_w - B_z -1\big) \eeq
For the component $\epsilon_{wx}$ this is immediate. The expression for $\epsilon_{yx}$ holds since $\beta_{xy}=\beta_{xu} -1 $ modulo $\alpha_y$, and the expression for $\epsilon_{wz}$ is derived similarly. The remaining three components also involve these kind of relations but all the $\pm 1$'s which appear cancel out.

Now simply observe that \eqref{eq.epsilon} gives a well-defined element:
$$\hat{\epsilon} \;\in \frac{ \Z_4^{\oplus 6}}{\smat{ 0 & 2 &2& 2 & 2 & 0\\ 2 & 0 & 2 & 2 & 0 & 2  }}  $$
We define $\hat{\epsilon}'$ in the same way, and it's immediate that if $\beta$ and $\beta'$ are equivalent then $\hat{\epsilon}= \hat{\epsilon}'$. 
For later use we record the following two identities that hold between the components of $\hat{\epsilon}$:

\beq{eq.refeps1}   \hat{\epsilon}_{wz} - \hat{\epsilon}_{yx} =   \hat{\epsilon}_{wy}- \hat{\epsilon}_{zx} + 2 \quad\; \in\Z_4 \eeq

\beq{eq.refeps2}  \hat{\epsilon}_{yx} + \hat{\epsilon}_{wx} = \hat{\epsilon}_{zx} + \hat{\epsilon}_{wy} - 2 \hat{\epsilon}_{zy} \quad \in \Z_4 \eeq

Note that the individual components of $\hat{\epsilon}$ only make sense mod 2, but the sums appearing in these identities all make sense mod 4. 
\end{rem}

\subsection{The refined flop equivalence}

\begin{prop}\label{prop.delicate} Let $\Gamma$ and $\Gamma'$ be the two graphs in Figure \ref{fig.IHmoves}. Let $(\alpha, \beta)$ be a decoration on $\Gamma$, and let $(\alpha, \beta')$ be an equivalent decoration on $\Gamma'$. Then the categories $\cC_\Gamma$ and $\cC_\Gamma'$ are equivalent relative to their external edges. 
\end{prop}

This is a more refined version of Proposition \ref{prop.5dflop2}, and to prove it we just need to re-examine the constructions in Section \ref{sec.flops} more carefully.  We recall all our notation from that section, in particular the 5-fold flop $Y \dashrightarrow Y'$. We fix a R-charge on $Y$ and $Y'$ matching our choice of $\alpha$. 

Suppose we have chosen an equivalence
$$K : D^b_\alpha(Y, W) \isoto \cC_\Gamma$$
as in Lemma \ref{lem.5dflop1}. This is the data of three Kn\"orrer equivalences
$$K_1: D^b_\alpha(Y_{s\neq 0}, W) \isoto D^b_\alpha(\A^3, uxy), \hspace{2cm} K_2: D^b_\alpha(Y_{t\neq 0}, W) \isoto D^b_\alpha(\A^3, vzw) $$
 and:
$$K_3: D^b_\alpha(Y_{st\neq 0}, W) \isoto D^b_\alpha(\C[u^{\pm 1}]) $$
Each of these equivalences on their own are unique up to shifts. To fit together to form the equivalence $K$ they must commute with the functors 
$$D^b_\alpha(\A^3, uxy) \to D^b_\alpha(\C[u^{\pm 1}])\aand D^b_\alpha(\A^3, vzw) \to D^b_\alpha(\C[u^{\pm 1}])$$
 specified by the decoration $\beta$, this means that $K_3$ is determined by either $K_1$ or $K_2$, and implies a certain compatibility between $K_1$ and $K_2$. For example, consider the object
$$\cO_Y/(y,z) \in D^b_\alpha(Y, W)$$
\emph{i.e.} the sky-scraper sheaf along the subvariety $\{y=z=0\}\subset Y$. Note that this is indeed an object in $D^b_\alpha(Y,W)$ because $W$ vanishes along this subvariety. We can apply $K_1$ or $K_2$ to (the restrictions of) this sheaf, and we must have
$$K_1(\cO_Y/(y,z)) = \cO/y[k_1] \aand K_2(\cO_Y/(y,z)) = \cO/z[k_2] $$ 
for some integers $k_1$ and $k_2$. For the diagram to commute we must have that:
\begin{align*} &&k_1  + \beta_{ux} &\equiv k_2 + \beta_{vw} \mod \alpha_u\\
 &\iff & k_2 - k_1 &\equiv \delta_{xw}  \end{align*}
Moreover, any pair $(k_1, k_2)$ obeying this congruence specifies an equivalence $K$. 

Now let us show how we can detect our element $B(\beta)$. Our decoration $\beta$ also specifies functors 
$$D^b_\alpha(\A^3, uxy) \to D^b_\alpha(\C[x^{\pm 1}])\aand D^b_\alpha(\A^3, vzw) \to D^b_\alpha(\C[w^{\pm 1}])$$
corresponding to the external edges $x$ and $w$ of $\Gamma$. If we pass our sheaf $\cO_Y/(y,z)$ through $K$ and then through these functors we get:
$$ \cO[k_1 + \beta_{xu}] \in  D^b_\alpha(\C[x^{\pm 1}])\aand \cO[k_2 + \beta_{wv}]\in D^b_\alpha(\C[w^{\pm 1}])$$
These two residues will detect two components of $B(\beta)$. But at the edges $y$ and $z$ the sheaf $\cO/(y,z)$ will go to zero, so to detect the remaining two components we need to consider at least one more object of $D^b_\alpha(Y,W)$. It is convenient to consider two more sky-scraper sheaves:
$$\cO_Y/(t, z) \aand \cO_Y/(y, s) \; \; \in D^b_\alpha(Y, W)$$
They are supported inside the two toric charts $\{s\neq 0\}$ and $\{t\neq 0\}$ respectively, and we must have:
$$K_1(\cO_Y/(t,z)) = \cO/u[k_1] \aand K_2(\cO_Y/(y,s) = \cO/v[k_2] $$
Using the functors at the external edges $y$ and $z$ specified by $\beta$ then gives:
$$ \cO[k_1+ \beta_{yx}]  \in  D^b_\alpha(\C[y^{\pm 1}])\aand \cO[k_2 + \beta_{zw}]\in D^b_\alpha(\C[z^{\pm 1}])$$
Our three sky-scraper sheaves have determined four residues
$$\big( k_1 + \beta_{xu},\,k_1 + \beta_{yx}, \, k_2 + \beta_{zw}, \, k_2 + \beta_{wv}\big)\;\in \Z_{\alpha_x}\oplus \Z_{\alpha_y}\oplus\Z_{\alpha_z}\oplus\Z_{\alpha_w} $$
and if consider the class of this element in the quotient $G_\alpha$ we get exactly $B(\beta)$. 

\begin{rem}\label{rem.lbonY} To pass to $G_\alpha$ we quotient by the elements $(1,1,1,1)$ and $(0,0,\alpha_u, \alpha_u)$. The first corresponds the shift autoequivalence of $D^b_\alpha(Y, W)$, since if we compose $K$ with a shift it adds 1 to all four residues. The second corresponds to tensoring by the `line-bundle' $O(1)$ on $Y$. This is an object which is $\cO$ on the chart $\{s\neq 0\}$ and $\cO[\alpha_u]$ on the chart $\{t\neq 0\}$, glued by the transition function $u=t/s$ (\emph{c.f.}~Remark \ref{rem.bmodalpha}).
 \end{rem}

\begin{lem} The equivalence $\Phi: D^b_\alpha(Y, W) \isoto D^b_\alpha(Y', W) $ from Theorem \ref{thm.VGIT} can be chosen such that:
$$ \cO_Y/(y,z) \mapsto \cO_{Y'}/(y,z), \quad\;\; \cO_Y/(t,z) \mapsto \cO_{Y'}/(t,z), \quad\;\;\cO_Y/(y,s) \mapsto \cO_{Y'}/(y,s) $$
\end{lem}
\begin{proof}This follows quite easily from the construction of $\Phi$ in \cite{segal} using `windows'. Recall that $Y$ and $Y'$ are quotients of $\A^6$ with co-ordinates $(s,t,p,q, y,z)$ by a $\C^*$ action of weights $(1,1,-1,-1,0,0)$, and that $W= tpy+sqz$.  

The functor $\Phi$ is induced from a derived equivalence $D^b_\alpha(Y)\isoto D^b_\alpha(Y')$ that sends $\cO\mapsto \cO(w)$ and $\cO(-1)\mapsto \cO(1+w)$, for some $w\in \Z$. We choose the $\Phi$ with $w=0$. 

The object $\cO_Y/(t, z) \in D^b_\alpha(Y, W)$ is equivalent to the following `perturbed Kozsul resolution':
$$\begin{tikzcd} \cO(-1) \arrow[yshift=0.7ex]{r}{(-z, t)} \ & \cO(-1)\oplus \cO \arrow[yshift=0.7ex]{r}{\smat{t \\ z}}  
\arrow[yshift=-0.7ex]{l}{\smat{-sq \\ py}}
& \cO \arrow[yshift=-0.7ex]{l}{(py, sq)} 
\end{tikzcd}$$
Applying $\Phi$ produces the same matrix factorization on $Y'$ but with $\cO(-1)$ replaced by $\cO(1)$. Hence $\Phi: \cO_Y/(t,z))\mapsto \cO_{Y'}/(t,z)$.

The same argument works for the other two objects. Note that the resolution for $\cO_Y/(y,z)$ only uses $\cO$ and not $\cO(-1)$.
\end{proof}

On the $Y'$ side, any two Kn\"orrer equivalences $K_1'$ and $K_2'$ on the two toric charts must obey
$$K_1'(\cO_{Y'}/(y,z)) = \cO[k_1'], \quad\quad K_2'(\cO_{Y'}/(y,z)) = \cO[k_2'] $$
and they patch to give an equivalence $K': D^b_\alpha(Y', W) \isoto \cC_{\Gamma'}$ iff:
$$k_2' - k_1' \equiv \delta'_{xw} \mod \alpha_{u'} $$
Applying $K'$ to the other two sheaves $\cO_{Y'}/(t,z)$ and $\cO_{Y'}/(y,s)$ lets us detect 
$$\big(k_1' + \beta'_{xu'},\, k_2'+\beta'_{yv'},\, k_1'+\beta_{zu'},\, k_2'+\beta'_{wv'}\big) \;\in \Z_{\alpha_x}\oplus \Z_{\alpha_y}\oplus\Z_{\alpha_z}\oplus\Z_{\alpha_w} $$
and the class of this element in $G'_\alpha$ is $B'(\beta')$. 

\begin{proof}[Proof of Proposition \ref{prop.delicate}] 
For any choice of $K$ and $K'$, the composition $\Psi = K'\circ \Phi\circ K^{-1}$ is an equivalence from $\cC_\Gamma$ to $\cC_{\Gamma'}$. We already know (Proposition \ref{prop.5dflop2}) that $\Psi$ commutes, up to some shifts, with the restriction functors at external edges. It remains to show that we can choose $K$ and $K'$ such that the shifts are zero.

We can check the shifts by considering the three objects:
$$K^{-1}(\cO_Y/(y,z)), \quad K^{-1}(\cO_Y/(t, z)), \quad K^{-1}(\cO_Y/(y, s))$$
If $B(\beta)$ and $B'(\beta')$ define the same class in $\overline{G}_\alpha$ then we can lift them to the same element of $\Z_{\alpha_x}\oplus \Z_{\alpha_y}\oplus\Z_{\alpha_z}\oplus\Z_{\alpha_w}$, and from the preceding discussion this corresponds to a choice of $k_1, k_2, k_1', k_2'$ giving us our required $\Psi$.\end{proof}

\begin{rem}\label{rem.spherical}
On $Y'$ we have an autoequivalence given by tensoring with the line-bundle $\cO(1)$, as in Remark \ref{rem.lbonY}. If we conjugate this with $\Phi$ we get an autoequivalence of $D^b_\alpha(Y,W)$, it's a family spherical twist along the subvariety $\P^1\times \A^2 \subset Y$. This autoequivalence acts as the identity away from  $\P^1\times 0$ so it induces an autoequivalence of $\cC_\Gamma$ relative to its external edges. 

 This autoequivalence doesn't come from a trivial modification of the $\beta$ data, but it does explain the coarser equivalence relation from Remark \ref{rem.coarserER}.
 \end{rem}

\subsection{Genus zero or one}\label{sec.genus0or1}

If $\Sigma$ is a punctured surface of genus zero then a line field on $\Sigma$ is determined, up to the action of the mapping class group, by its winding numbers around the boundary curves \cite{LPinvs}. Therefore the Fukaya category only depends on these numbers. Our version of this statement is the following:

\begin{thm}\label{thm.genus0} Let $\Gamma$ and $\Gamma'$ be two connected genus zero trivalent graphs with a chosen bijection $f: \partial \Gamma \to \partial \Gamma'$. Choose decorations $(\alpha, \beta)$ and $(\alpha', \beta')$ on each graph such that $\alpha_x = \alpha'_{f(x)}$ for each external edge $x$. Then we have an equivalence:
$$\cC_\Gamma \isoto \cC_{\Gamma'}$$
\end{thm}


Note that fixing the values of $\alpha$ on $\partial \Gamma$ determines $\alpha$ completely.

\begin{proof}
By Lemma \ref{lem.IHmoves} there is a sequence of IH moves transforming $\Gamma$ into $\Gamma'$. By Proposition \ref{prop.5dflop2} after each IH move there is some decoration on the new graph such that the limit categories are unchanged. Hence there is some decoration $(\alpha', \beta'')$ on $\Gamma'$ for which $\cC_{\Gamma'}$ is equivalent to $\cC_{\Gamma}$. But by Lemma \ref{lem.genus0trivmod} the category $\cC_{\Gamma'}$ only depends on the $\alpha$ part of the decoration.
\end{proof}

Now let $\Sigma$ be a punctured surface of genus one, equipped with a line field $\eta$. Taking the winding numbers of $\eta$ around the boundary components gives integers $\alpha_1,..., \alpha_k$. Now take two more embedded curves which, after filling in the punctures, give a basis of $H^1$ of the torus; the winding numbers around these curves give two more integers $a, b$. The orbit of $\eta$ under the mapping class group is determined  \cite{LPinvs} by the $\alpha_i$ together with the invariant:
$$ \widetilde{A}(\eta) = \gcd( \alpha_1 +2, ..., \alpha_k + 2, a, b) $$
Now let $\Gamma$ be a connected trivalent graph of genus one with a decoration $(\alpha, \beta)$. Obviously we can get integers 
$$\alpha_1,..., \alpha_k \in \Z$$
by taking the values of $\alpha$ on each external edge, but we must also explain how to extract this invariant $\widetilde{A}$. 

The graph $\Gamma$ consists of a single cycle with trees attached to it. Let $\Delta$ be the set of half-edges which meet a vertex appearing in the cycle, so $\Delta$ forms a genus one graph with as many external edges as vertices (in Section \ref{sec.STZ} we called this a `wheel').  Now let $x$ be an external edge of $\Delta$,  it is one leaf of a tree whose remaining leaves $y,..., z$ lie in $\partial \Gamma$. It follows immediately from \eqref{eq.eulerchar} that:
$$\alpha_x - 2 = (\alpha_y-2) + ... + (\alpha_z-2)$$
Now let $u$ be half of an internal edge of $\Delta$, \emph{i.e.} a half-edge of $\Gamma$ that is part of the cycle. Obviously if $v$ a different such half-edge then $\alpha_u$ and $\alpha_v$ may not be equal, but the formula above implies that the residue 
$$a = [\alpha_u] = [\alpha_v]  \in \Z/(\alpha_1 -2, ..., \alpha_k - 2) $$
is well defined. 

  If $u,v$ are two such half-edges that meet at a vertex then we have a quantity $\gamma_{uv} = \beta_{uv} - \beta_{vu} \in \Z/(\alpha_u, \alpha_v)$, so summing these around the cycle as in \eqref{eq.bc} gives a well-defined residue:
 $$b \in  \Z/(\alpha_1 -2, ..., \alpha_k - 2, a) $$
Our required invariant is
$$\widetilde{A}(\alpha, \beta) = \gcd( \alpha_1-2, ..., \alpha_k - 2, a, b) $$
 (the switch to $\alpha_i-2$ from the $\alpha_i+2$ used in \cite{LPinvs} is just an orientation convention). Since $b$ is invariant under trivial modifications $\widetilde{A}$ is too. 

\begin{thm}\label{thm.genus1} Let $\Gamma$ and $\Gamma'$ be two connected genus one trivalent graphs with a chosen bijection $f: \partial \Gamma \to \partial \Gamma'$. Choose decorations $(\alpha, \beta)$ and $(\alpha', \beta')$ on each graph, such that $\alpha_x = \alpha'_{f(x)}$ for each external edge $x$, and such that $\widetilde{A}(\alpha, \beta) = \widetilde{A}(\alpha', \beta')$. Then we have an equivalence:
$$\cC_\Gamma \isoto \cC_{\Gamma'}$$
\end{thm}
\begin{proof}
Suppose we perform an IH move on $\Gamma$, and add a decoration to the new graph which is equivalent to $(\alpha, \beta)$. Neither $a$ nor $b$ change under this move (see Remark \ref{rem.bcunderIH}) hence $\widetilde{A}$ doesn't change either.

Now let $\Gamma''$ be a third graph which is a tree with a single loop attached, and with external edges $\partial \Gamma''=\partial \Gamma=\partial \Gamma'$.   By Lemma \ref{lem.IHmoves} we can apply IH moves to turn either $\Gamma$ or $\Gamma'$ into $\Gamma''$. Note that on $\Gamma''$ the $\widetilde{A}$ invariant is particularly simple: $a$ is just the value of $\alpha$ on the loop and makes sense as an integer, and $b$ makes sense in $\Z_a$.

Suppose we start at $\Gamma$ with the decoration $(\alpha, \beta)$, and after every IH move replace our decoration with an equivalent decoration. We will finish with some decoration on $\Gamma''$ having the same $\widetilde{A}$ value as we started with. By Proposition \ref{prop.delicate} this decoration produces a category equivalent to $\cC_\Gamma$. 

Alternatively we start on $\Gamma'$ with the decoration $(\alpha', \beta')$ and end with some other decoration on $\Gamma''$, giving a category equivalent to $\cC_{\Gamma'}$. These two decorations on $\Gamma''$ have the same $\widetilde{A}$ invariant, so it follows from Corollary \ref{cor.SL2} and Lemma \ref{lem.FMcommuteswithDsg} that they are equivalent.
\end{proof}

\subsection{Genus at least two}\label{sec.genus2invariance}

Let $\Gamma$ be a decorated trivalent graph. Let us suppose that
\begin{enumerate}
\item $\alpha_x$ is even for every half edge $x$ in $\Gamma$.
\item $\alpha_x\equiv 2$ mod 4 for every external edge $x\in \partial \Gamma$. 
\end{enumerate}
Then we may define 
$$ p_x = \frac12 \alpha_x + 1 \;\; \in \Z_2$$
for each half edge. Since these residues sum to zero at each vertex and vanish on external edges we have defined a class in $H_1(\Gamma, \Z_2)$. This is exactly the data of a set of pair-wise disjoint cycles $\mathfrak{C}$ in $\Gamma$. A cycle $c$ lies in $\mathfrak{C}$ iff $\alpha_x\equiv 0$ mod 4 for all half-edges $x\in c$. 

Now take any cycle $c$. From the $\beta$ data we get an invariant $b_c$ as in \eqref{eq.bc}, which lives in the quotient of $\Z$ by the ideal $I_c$ generated by $\alpha_x$ for each half-edge $x$ in the cycle. Let's further assume that
\begin{enumerate}\addtocounter{enumi}{2}
\item $b_c \equiv 0 $ mod 2 for all cycles $c$.
\end{enumerate}
which makes sense because of assumption (1). 
\begin{rem}\label{rem.associatedweak2} Another way to say this is that (under assumption (1)) we can consider the associated weak decoration we get by reducing $\beta$ mod 2. Then assumption (3) says that this weak decoration defines the zero class in $H^1(\Gamma, \Z_2)$. 
\end{rem}

If $c\in \mathfrak{C}$ then we may consider the value of $b_c$ mod 4, and then the quantity
$$ q_c = \frac12 b_c + 1\; \in \Z_2 $$
is well-defined (and independent of the orientation of $c$). We define:
$$A(\alpha, \beta) = \sum_{c\in \mathfrak{C}} q_c \quad \in \Z_2$$
\begin{rem} This agrees with the Arf invariant of a line field on a surface; see \cite[Sec.~1.2]{LPinvs} and the references therein.  Suppose we have a line field $\eta$ on a surface $\Sigma$, and a pair-of-pants decomposition, resulting in our decorated graph $\Gamma$. Extend $\mathfrak{C}$ to a basis of $H_1(\Gamma)$ and then extend this to a symplectic basis of $H_1(\overline{\Sigma})$, where $\overline{\Sigma}$ is the compact surface obtained by filling in the punctures. Taking winding numbers of $\eta$ along these cycles gives integers $(\tilde{\alpha}_1, \tilde{\beta}_1, ..., \tilde{\alpha}_g, \tilde{\beta}_g)$ which are all even from our assumptions on the decoration.  We set $p_i = \tfrac12 \tilde{\alpha}_i + 1 \in \Z_2$ and  $q_i = \tfrac12 \tilde{\beta}_i + 1$. The Arf invariant of the quadratic form
$$\sum_{i=1}^g p_iX_i^2 + X_iY_i + q_j Y_i^2 $$
is $A= \sum_{i=1}^g p_iq_i$.
 
The cycles in $\mathfrak{C}$ give $p_i=1$ and for these cycles $q_i$ agrees with our $q_c$. For the remaining cycles in $H_1(\Gamma)$ we can arrange that $p_i=0$ so these do not contribute to $A$.
\end{rem}

\begin{rem}\label{rem.assumptionsunderIH} Suppose $\Gamma$ and $\Gamma'$ are two graphs related by an IH move and $(\alpha, \beta)$, $(\alpha', \beta')$ are two equivalent decorations. If $(\alpha, \beta)$ satisfies assumptions (1)-(3) then so does $(\alpha', \beta')$. For (1) and (2) this is obvious and for (3) it follows from Remarks \ref{rem.associatedweak1} and \ref{rem.associatedweak2}.
\end{rem}

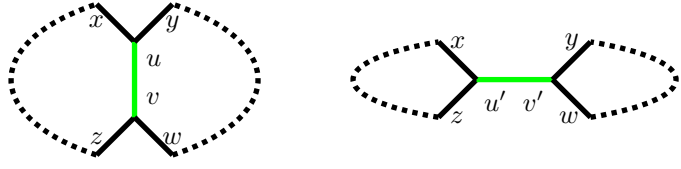
\begin{figure}\label{fig.IIIa}
\begin{tikzpicture}[scale =.5]
\draw (0,0) --(1,1);
\draw (1,3)--(0,4);
\draw[green] (1,1)--(1,3);
\draw (1,3)--(2, 4);
\draw (1,1)--(2,0);
\node [below] at (0,4) {$x$};
\node [below] at (2,4) {$y$};
\node [right] at  (1, 2.5) {$u$};
\node [right] at  (1, 1.5) {$v$};
\node [above] at  (0, 0) {$z$};
\node [above] at  (2,0) {$w$};
\draw[dotted] (0,0) .. controls (-3,1) and (-3, 3) .. (0,4);
\draw[dotted] (2,0) .. controls (5,1) and (5,3).. (2, 4);

\begin{scope}[shift={(4,0)}]
\draw (5,1) --(6,2);
\draw[green] (6,2)--(8,2);
\draw (8,2)--(9,1);
\draw (6,2)--(5, 3);
\draw (8,2)--(9,3);
\node [right] at (5,3) {$x$};
\node [left] at (9,3) {$y$};
\node [below] at  (6.5, 2) {$u'$};
\node [below] at  (7.5, 2) {$v'$};
\node [right] at  (5, 1) {$z$};
\node [left] at  (9,1) {$w$};
\draw[dotted] (5,1) .. controls (2,1.5) and (2,2.5) .. (5,3);
\draw[dotted] (9,1) .. controls (12,1.5) and (12,2.5).. (9, 3);
 
  \end{scope}
\end{tikzpicture}
\caption{Subcase (i)}\label{fig.subcase1}
\end{figure}

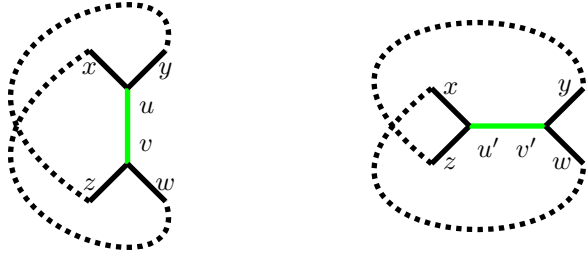
\begin{figure}\label{fig.IIIb}
\begin{tikzpicture}[scale =.5]
\draw (0,0) --(1,1);
\draw (1,3)--(0,4);
\draw[green] (1,1)--(1,3);
\draw (1,3)--(2, 4);
\draw (1,1)--(2,0);
\node [below] at (0,4) {$x$};
\node [below] at (2,4) {$y$};
\node [right] at  (1, 2.5) {$u$};
\node [right] at  (1, 1.5) {$v$};
\node [above] at  (0, 0) {$z$};
\node [above] at  (2,0) {$w$};
\draw[dotted] (0,0) .. controls (-6,4) and (3,7)   .. (2,4);
\draw[dotted] (0,4) .. controls (-6,0) and (3,-3).. (2, 0);

\begin{scope}[shift={(4,0)}]
\draw (5,1) --(6,2);
\draw[green] (6,2)--(8,2);
\draw (8,2)--(9,1);
\draw (6,2)--(5, 3);
\draw (8,2)--(9,3);
\node [right] at (5,3) {$x$};
\node [left] at (9,3) {$y$};
\node [below] at  (6.5, 2) {$u'$};
\node [below] at  (7.5, 2) {$v'$};
\node [right] at  (5, 1) {$z$};
\node [left] at  (9,1) {$w$};
\draw[dotted] (5,1) .. controls (0,5) and (9,6) .. (9,3);
\draw[dotted] (5,3) .. controls (0,-1) and (9,-2).. (9, 1);
 
  \end{scope}
\end{tikzpicture}
\caption{Subcase (ii)}\label{fig.subcase2}
\end{figure}

\begin{prop}\label{prop.AisIHinvariant} Let $\Gamma$ and $\Gamma'$ be two graphs related by an IH move. Let $(\alpha, \beta)$ be a decoration on $\Gamma$ and let $(\alpha', \beta')$ be an equivalent decoration on $\Gamma'$. Assume that $(\alpha, \beta)$ - and hence $(\alpha', \beta')$ - satisfy the conditions (1)-(3) above. Then:
$$ A(\alpha, \beta) = A(\alpha', \beta') $$
\end{prop}
\begin{proof}
As usual consider the half-edges in $\Gamma$ and $\Gamma'$ involved in the IH move and label them as in Figure \ref{fig.IHmoves}. Since $\alpha_x+\alpha_y+\alpha_z + \alpha_w=4$ we must have that either 0, 2 or 4 of these values are congruent to zero mod 4. We deal with each case separately, splitting the third case into two sub-cases.
\pgap

\noindent \emph{None are congruent to zero mod 4.} In this case no cycle in $\mathfrak{C}$ is affected by the IH move and $A$ doesn't change.
\pgap

\noindent  \emph{Two are congruent to zero mod 4.}  Then at least one of $\alpha_u$ and $\alpha'_{u'}$ equals 0 mod 4. In this case there is one cycle $c\in \mathfrak{C}$ affected by the IH move, and it is transformed to a cycle $c'\subset \Gamma'$ which lies in the corresponding set $\mathfrak{C}'$. At least one of $c$ or $c'$ includes the edge affected by the IH move. By Remark \ref{rem.bcunderIH} we have $b_c \equiv b_{c'}$ mod 4, hence $A$ does not change. 
\pgap

\noindent \emph{All four values are zero mod 4.} Then $\alpha_u\equiv \alpha'_{u'}\equiv 2$ mod 4. Either we have two cycles $c_1, c_2\in \mathfrak{C}$ containing $x, y$ and $z,w$ respectively, or we have a single cycle $c\in \mathfrak{C}$ containing $x,y,z, w$. In the first case the IH move merges the two cycles to a single cycle $c'\in \mathfrak{C}$ such that
$$ c' = c_1 + c_2 \; \in H_1(\Gamma, \Z_2) $$
(using the canonical isomorphism between $H^1(\Gamma)$ and $H^1(\Gamma')$). In the second case the IH move may split $c$ into $c_1'+c_2'$, or it may turn it into a single cycle $c'$. So by symmetry the two subcases we need to deal with are $c\leadsto c'_1+c'_2$ or $c\leadsto c'$. See Figures \ref{fig.subcase1} and \ref{fig.subcase2}. 
\pgap

\noindent \emph{Subcase (i).} $c\leadsto c'_1 + c'_2.\quad$ The quantity $b_c$ can be expressed as
$$b_c = \gamma_{xy} + b_{yw} + \gamma_{wz} + b_{zx} $$
where $b_{yw}$ and $b_{zx}$ denote the sum of the $\gamma$ invariants along the two segments of $c$ which are not affected by the IH move. After the move we have:
$$b_{c_1'} = \gamma'_{xz} + b_{zx} \aand b_{c_2'} = \gamma'_{wy} + b_{yw} $$
We are interested in the residue of these quantities in $\Z_4$. Using our refined invariants $\hat{\epsilon}$ and $\hat{\epsilon}'$ from Remark \ref{rem.refinedepsilon} we have:
$$ b_c = \hat{\epsilon}_{wz}-\hat{\epsilon}_{yx} + b_{yw}+ b_{zx} \; \in \Z_4$$
and
$$b_{c_1'} + b_{c_2}' = \hat{\epsilon}'_{wy} - \hat{\epsilon}'_{zx} + b_{zx} + b_{yw}\; \in \Z_4$$
Since $\beta$ and $\beta'$ are equivalent we have $\hat{\epsilon}= \hat{\epsilon}'$, and now the identity \eqref{eq.refeps1} implies that:
$$b_c = b_{c_1'} + b_{c_2'} + 2 \quad \in \Z_4 $$
It follows that $q_c = q_{c_1'} + q_{c_2'}$, and $A$ does not change.
\pgap

\noindent \emph{Subcase (ii).} $c\leadsto c'. \quad$ As above the  quantities $b_c$ and $b_{c'}$ can be expressed as
$$ b_c = \gamma_{yx} + b_{xw} + \gamma_{wz} + b_{zy} = \hat{\epsilon}_{yx} + \hat{\epsilon}_{wz} + b_{xw}- b_{yz} \; \in \Z_4 $$
and
$$ b_{c'} = \gamma'_{zx} + b_{xw} + \gamma'_{wy} - b_{zy} = \hat{\epsilon}'_{zx} + \hat{\epsilon}'_{wy} + b_{xw }+b_{yz} \; \in \Z_4 $$
Again $\hat{\epsilon}= \hat{\epsilon}'$, and now the identity \eqref{eq.refeps2} implies that
$$b_{c} = b_{c'} - 2\hat{\epsilon}_{zy}  - 2 b_{yz}\; \in \Z_4 $$
The quantity $\hat{\epsilon}_{zy}+ b_{yz} \in \Z_2$ is an invariant computed around a cycle so must be zero by assumption (3). So $b_c= b_{c'}\in \Z_4$ and hence $A$ does not change. 
\end{proof}

We are almost ready to state the analogue of Theorems \ref{thm.genus0} and \ref{thm.genus1}. For a graph $\Gamma$ we split the set of all decorations on $\Gamma$ into the following four disjoint classes, corresponding to the assumptions (1)-(3):
\begin{itemize}\item[(I)] Either the values of $\alpha$ include at least one odd number, or all  $\alpha$ values are even and $b_c$ is odd for at least one cycle.
\item[(II)] Every value of $\alpha$ and every $b_c$ is even. At least one external edge $x\in \partial\Gamma$ satisfies $\alpha_x\equiv 0$ mod 4. 
\item[(III)] Assumptions (1)-(3) hold and $A(\alpha, \beta)=0$.
\item[(IV)] Assumptions (1)-(3) hold and $A(\alpha, \beta)=1$.
\end{itemize}

Each of these classes is preserved if we perform an IH move on $\Gamma$ and equip the new graph with an equivalent decoration; see Remark \ref{rem.assumptionsunderIH} and Proposition \ref{prop.AisIHinvariant}.

\begin{thm} Let $\Gamma$ and $\Gamma'$ be two connected trivalent graphs of genus $g\geq 2$ with a chosen bijection $f: \partial \Gamma \to \partial \Gamma'$. Choose decorations $(\alpha, \beta)$ and $(\alpha', \beta')$ on each graph such that $\alpha_x = \alpha'_{f(x)}$ for each external edge $x$. If $(\alpha, \beta)$ and $(\alpha', \beta')$ lie in the same class $I-IV$ then we have an equivalence:
$$\cC_\Gamma \isoto \cC_{\Gamma'}$$
\end{thm}

This is similar to \cite[Thm.~1.2.4(iii)]{LPinvs} and the earlier results which that is modelled on. We give a complete proof here since it is elementary and fairly short.

\begin{proof} The first steps are similar to our proof of Theorem \ref{thm.genus1}. By performing IH moves we may assume (without changing either category) that $\Gamma$ and $\Gamma'$ are the same graph, and we may choose this graph to be a tree with $g$ loops attached. The values of $\alpha$ on the tree are determined by the values on $\partial \Gamma$, the remaining choices are the values $\alpha_1,..., \alpha_g$ on each loop. The $\beta$ data gives us an invariant $b_i\in \Z_{\alpha_i}$ for each loop. Choose a lift $(\tilde{b}_1,..., \tilde{b}_g)\in \Z^g$. The category $\cC_\Gamma$ is determined by these $2g$ integers.

Similarly the category $\cC_{\Gamma'}$ is determined by $2g$ integers $(\alpha'_1,..., \alpha'_g)$ and $(\tilde{b}'_1,..., \tilde{b}'_g)$. We have various operations that we can perform on these tuples of integers without affecting the categories; our task is to show that together these operations act transitively within each class I-IV. 

The obvious operations are that we may change $\tilde{b}_i$ by a multiple of $\alpha_i$ for each $i$, we can multiply any pair $(\alpha_i, \tilde{b}_i)$ by $-1$, and we may perform any permution of the set $[1,g]$. We need to consider four more operations:
\begin{itemize}\setlength{\itemsep}{10pt}

\item For each $i$ we may change $(\alpha_i, \tilde{b}_i)$ into $(\tilde{b}_i, -\alpha_i)$. This doesn't change the category by Proposition \ref{prop.FM} and Lemma \ref{lem.FMcommuteswithDsg}.

\item We may change the value of $\alpha_1$ (hence any $\alpha_i$) by a multiple of 4. To see that this doesn't affect the category inspect the graph on the LHS of Figure \ref{fig.doubleapple2}; we can perform IH moves until this appears as a subgraph, then we can swap $\alpha_1$ with $\alpha_1-4$.  

Obviously $\alpha_i, \tilde{b}_i$ are not affected for $i>1$. The invariant $b_1\in \Z/\alpha_1$ changes to an invariant $b_1'\in \Z/(\alpha_1-4)$, but since these must be the same in $\Z/(\alpha_1, 4)$ we can lift them to the same integer $\tilde{b}_1$. 

\item For any external edge $x$ we may change the value of $\alpha_1$ (hence any $\alpha_i$) by a multiple of $\alpha_x-2$. This is similar to the previous operation, perform IH moves until the middle of Figure \ref{fig.doubleapple1} appears. 

\item We can change $(\alpha_1,\tilde{b}_1, \alpha_2, \tilde{b}_2)$ to:
$$(\alpha_1, \tilde{b}_1 + \tilde{b}_2 - 2, \alpha_2 - \alpha_1 + 2, \tilde{b}_2)$$
See Figure \ref{fig.doubleapple2}.  To get to the right-hand-side we have first perfomed an IH move, then changed basis on the cycles from $(c_1, c_2)$ to $(c_1+ c_2, c_2)$. The change in $\tilde{b}_1$ is explained by the fact that $b_{c_1+c_2} \equiv b_{c_1} + b_{c_2} - 2 $ 
\end{itemize}

Now we can deal with each class of decoration.

\begin{itemize}\setlength{\itemsep}{10pt}
\item[\emph{Class (I):}] If $\alpha_1$ (say) is odd we can reduce it to $\pm 1$, then we can reduce all the $\alpha_i$ to 1. Then the $\tilde{b}_i$ are irrelevant. If all $\alpha_i$ are even but $\alpha$ is odd on some half-edge then $\alpha_x$ must be odd for at least one external edge. Then we can make $\alpha_1$ odd. 
 
 Now suppose all $\alpha$ values are even but $b_c$ is odd for some cycle. From the proof of Lemma \ref{lem.IHmoves} we can arrange (by cutting the cycle $c$ first) that $b_1$ is odd. Then we can swap it with $\alpha_1$.

\item[\emph{Class (II):}] We can reduce each $\alpha_i$ to 2 and then each $\beta_i$ to zero.

\item[\emph{Classes (III)+(IV):}] If $\alpha_i\equiv \tilde{b}_i\equiv 0 \mod 4$ then we can reduce this pair to $(0,0)$. If not then we can  reduce this pair to $(2,0)$. The number of each of kind of pairs is only relevant modulo 2, using the final operation listed above. So we can reduce to the case that $\tilde{b}_i =0$ for all $i$ and $\alpha_j=2$ for all $j>1$, and $\alpha_1$ equals either zero or 2 depending on $A(\alpha, \beta)$. 
\end{itemize}

\end{proof}

\begin{figure}
\begin{tikzpicture}[scale =.5]
\begin{scope}[decoration={
    markings,
    mark=at position 0.5 with {\arrow{>}}}
    ] 
\draw[postaction={decorate}] (1,0) --(-1,0);
\draw[postaction={decorate}] (5,0) --(3,0);
\draw[postaction={decorate}]  (3,0) arc [radius=1, start angle=0, end angle= 180];
\draw[postaction={decorate}]  (1,0) arc [radius=1, start angle=180, end angle= 360];
\draw[postaction={decorate}]  (5,0) arc [radius=1, start angle=-180, end angle= 180];
\node [below] at (-1,-.1) {6};
\node [below] at (4,-.1) {2};
\node [above] at (2,1.1) {$\alpha_1$};
\node[below] at (2, -1.1) {$\alpha_1-4$};
\node [right] at (7.1,0) {$\alpha_2$};
\end{scope}

\begin{scope}[shift={(13,0)},  decoration={
    markings,
    mark=at position 0.5 with {\arrow{>}}}
    ] 
\draw[postaction={decorate}] (1,0) --(-1,0);
\draw[postaction={decorate}] (3,1.8) --(1,0);
\draw[postaction={decorate}] (1,0) --(3,-1.8);
\draw[postaction={decorate}] (3,-1.8) --(3,1.8);
\draw[postaction={decorate}]  (3,1.8) arc [radius=1.8, start angle=90, end angle= -90];
\node [below] at (-1,-.1) {6};
\node [right] at (5,0) {$\alpha_2-\alpha_1+2$};
\node [right] at (3.1,0) {$\alpha_2$};
\node [above] at (1.6,1.1) {$\alpha_1$};
\node[below] at (1.1, -1) {$\alpha_1-4$};
\end{scope}

\end{tikzpicture}
\caption{}\label{fig.doubleapple2}
\end{figure}


\begin{figure}
\begin{tikzpicture}[scale =.5]
\begin{scope}[decoration={
    markings,
    mark=at position 0.5 with {\arrow{>}}}
    ] 
\draw[postaction={decorate}] (1,0) --(-1,0);
\draw[postaction={decorate}] (3,0) --(1,0);
\draw[postaction={decorate}] (1,2) --(1,0);
\draw[postaction={decorate}] (1,2) arc [radius=1, start angle=-90, end angle= 270];

\node [below] at (-1,-0.1) {$\alpha_x$};
\node [below] at (3,-.1) {$\alpha_x-4$};
\node [left] at (.9,1) {2};
\node [above] at (1,4.1) {$\alpha_1$};
\end{scope}

\begin{scope}[shift={(9,1.5)},  decoration={
    markings,
    mark=at position 0.5 with {\arrow{>}}}
    ] 
\draw[postaction={decorate}] (1,0) --(-1,0);
\draw[postaction={decorate}] (5,0) --(3,0);
\draw[postaction={decorate}]  (3,0) arc [radius=1, start angle=0, end angle= 180];
\draw[postaction={decorate}]  (1,0) arc [radius=1, start angle=180, end angle= 360];
\node [below] at (-1,-.1) {$\alpha_x$};
\node [below] at (5,-.1) {$\alpha_x-4$};
\node [above] at (2,1.1) {$\alpha_1$};
\node[below] at (2, -1.1) {$\alpha_1-\alpha_x+2$};
\end{scope}

\begin{scope}[shift={(20,3)}, decoration={
    markings,
    mark=at position 0.5 with {\arrow{>}}}
    ] 
\draw[postaction={decorate}] (1,0) --(-1,0);
\draw[postaction={decorate}] (3,0) --(1,0);
\draw[postaction={decorate}] (1,-2) --(1,0);
\draw[postaction={decorate}] (1,-2) arc [radius=1, start angle=90, end angle= 450];

\node [below] at (-1,-.1) {$\alpha_x$};
\node [below] at (3,-.1) {$\alpha_x-4$};
\node [left] at (.9,-1) {2};
\node [below] at (1,-4.1) {$\alpha_1-\alpha_x+ 2$};
\end{scope}

\end{tikzpicture}
\caption{}\label{fig.doubleapple1}
\end{figure}

\bibliographystyle{halphanum}

\end{document}